\documentclass{article}

\usepackage{amsmath, amssymb}
\usepackage[amsthm, thref, hyperref,thmmarks, amsmath]{ntheorem}
\usepackage{mathtools} 
\usepackage{accents}
\usepackage{stmaryrd}

\usepackage{geometry}
\usepackage{hhline}
\usepackage{float}
\usepackage{multicol}
\usepackage{booktabs}

\usepackage{graphicx}
\usepackage{tikz}
\usetikzlibrary{calc, shapes, plotmarks, patterns, positioning, decorations.pathreplacing, calligraphy, arrows.meta}
\usepackage{wrapfig}

\usepackage{cancel}
\usepackage[labelfont=bf]{subcaption}

\usepackage{tabularx}
\usepackage{array}
\usepackage{calc}

\usepackage{hyperref}

\newtheorem{Lemma}{Lemma}[section]
\newtheorem{Theorem}{Theorem}[section]

\theoremstyle{definition}
\newtheorem{Definition}{Definition}[section]
\newtheorem{Assumption}{Assumption}[section]

\newtheorem*{Remark}{Remark}

\theoremsymbol{\hfill\ensuremath{\pmb\square}}
\newtheorem*{Proof}{Proof}

\numberwithin{equation}{section}

\newcommand{\RR}{\mathbb{R}}

\newcommand{\VV}{\mathbb{V}}

\newcommand{\Mc}{\mathcal{M}}

\newcommand{\Hc}{\mathcal{H}}
\newcommand{\Ic}{\mathcal{I}}
\newcommand{\Rc}{\mathcal{R}}

\newcommand{\Xc}{\mathcal{X}}
\newcommand{\Zc}{\mathcal{Z}}

\newcommand{\wh}[1]{\widehat #1}
\newcommand{\x}{\mathbf{x}}

\definecolor{royalblue}{rgb}{0.25, 0.41, 0.88}

\newcommand*\diff{\mathop{}\!\mathrm{d}}

\newcommand*\Div{\mbox{div}}
\newcommand*\Curl{\mbox{curl}}

\newcommand{\cl}[1]{\overline{#1}}

\newcommand{\uu}{\underline{u}}
\newcommand{\uv}{\underline{v}}
\newcommand{\uw}{\underline{w}}

\newcommand{\uf}{\underline{f}}
\newcommand{\ug}{\underline{g}}
\newcommand{\ud}{\underline{d}}

\newcommand{\ulam}{\underline{\lambda}}
\newcommand{\umu}{\underline{\mu}}

\newcommand{\Ir}{\mathrm{I}}

\definecolor{myblue}{RGB}{30,90,160}
\hypersetup{
  colorlinks=true,        
  urlcolor=myblue,        
  citecolor=black,
  linkcolor=black
}

\makeatletter
\renewcommand{\@fnsymbol}{\@arabic}
\makeatother

\title{Dual-primal Isogeometric Tearing and Interconnecting Solvers for adaptively refined multi-patch configurations}
\author{}
\date{\today}

\begin{document}
\maketitle
\vspace{-2em}
\begin{center}
\begin{minipage}[t]{0.4\textwidth}
    \centering
    \textbf{Stefan Takacs}\\
    \footnotesize \href{mailto:stefan.takacs@jku.at}{stefan.takacs@jku.at}\\
    Institute for Numerical Mathematics,\\
    Johannes Kepler University (JKU)\\
    Altenberger Straße 69, 4040 Linz
\end{minipage}
\hfill
\begin{minipage}[t]{0.4\textwidth}
    \centering
    \textbf{Stefan Tyoler}\\
    \footnotesize \href{mailto:stefan.tyoler@ricam.oeaw.ac.at}{stefan.tyoler@ricam.oeaw.ac.at}\\
    Johann Radon Institute for Computational and Applied Mathematics (RICAM)\\
    Altenberger Straße 69, 4040 Linz
\end{minipage}
\end{center}

\begin{abstract}
Isogeometric Analysis is a variant of the finite element method, where spline functions are used for the representation of both the geometry and the solution. Splines, particularly those with higher degree, achieve their full approximation power only if the solution is sufficiently regular. Since solutions are usually not regular everywhere, adaptive refinement is essential. Recently, a multi-patch-based adaptive refinement strategy based on recursive patch splitting has been proposed, which naturally generates hierarchical, non-matching multi-patch configurations with T-junctions, but preserves the tensor-product structure within each patch.

In this work, we investigate the application of the dual-primal Isogeometric Tearing and Interconnecting method (IETI-DP) to such adaptive multi-patch geometries. We provide sufficient conditions for the solvability of the local problems and propose a preconditioner for the overall iterative solver. We establish a condition number bound that coincides with the bound previously shown for the fully matching case. Numerical experiments confirm the theoretical findings and demonstrate the efficiency of the proposed approach in adaptive refinement scenarios.
\end{abstract}

\section{Introduction}

Isogeometric Analysis (IgA) was introduced by Tom Hughes et al \cite{Hughes2005} with the goal of unifying Computer-Aided Design (CAD) and Finite Element Analysis through the use of spline-based basis functions, such as B-splines and Non-Uniform Rational B-Splines (NURBS), both for the representation of the geometry and the solution. Since only trivial geometries can be parameterized with a single smooth geometry function, the overall computational domain is usually represented as the union of non-overlapping, individually parameterized patches (multi-patch IgA). In IgA, spline degrees and spline smoothness can be freely chosen, which may be beneficial in several applications such as structural mechanics, fluid dynamics, and wave propagation, where high-order approximations and continuity are often advantageous.

High-order approximations can reveal their full approximation power only if the solution is smooth enough. Differential equations resulting from practical applications exhibit localized singularities or sharp gradients in the solution, often arising from geometric features such as re-entrant corners, boundary layers, or discontinuities in material coefficients, see, for example, \cite{Schwab1998}. For high-order approximations, this is even more pronounced and highlights the necessity of adaptive mesh refinement strategies, which selectively enhance the resolution in regions where the solution exhibits reduced regularity, while keeping the overall number of degrees of freedom manageable.
To efficiently steer the refinement process, we consider residual-based a posteriori error estimators \cite{Verfuerth2013} in combination with the D\"orfler marking strategy \cite{Doerfler1996, Praetorius2020}.

In standard IgA, tensor-product spline functions are used as basis functions for the parameter domain. The basis functions for the physical domain are then defined via the pull-back principle.
In order to facilitate adaptivity, unstructured spline constructions have been in the focus for years, which leads to a zoo of spline constructions, including Hierarchical B-splines (HB) and Truncated Hierarchical B-splines (THB), see \cite{Kraft1997,Giannelli2011, Giannelli2012}, (analysis-suitable) T-splines, see \cite{Dörfel2010, Sederberg2012}, and Locally Refined splines (LR), see \cite{Dokken2013, Dokken2014}. All of these spline variants allow localized refinement, while preserving properties such as linear independence, partition of unity, and local support. For all of these approaches, effort has to be put into an efficient implementation, particularly in the context of multi-patch geometries.

Recently, the authors proposed a multi-patch-based adaptive scheme \cite{Tyoler2025}, based on the repetitive splitting of existing patches into sub-patches in order to enable local refinement. Unlike other spline constructions, where extending to multi-patch domains often involves various compatibility conditions across interfaces, the application of this patch-splitting strategy to multi-patch geometries is straightforward by construction. The construction preserves the local tensor-product structure for each patch. Certainly, this comes at a price since the resulting discretization is not fully matching. However, the construction guarantees that the spaces at the interfaces are nested, which means that on each interface, the trace space of the spline space of one of the two adjacent patches is contained in the trace space of the spline space of the other adjacent patch. The interfaces between patches may result in T-junctions.

Efficient solvers are crucial for addressing the large, ill-conditioned linear systems arising from high-order isogeometric discretizations. While direct solvers remain a robust and widely used approach, especially for moderate problem sizes, their computational complexity and memory requirements render them less feasible in large-scale or three-dimensional settings. 
In these settings, iterative solvers become indispensable for achieving computational scalability. 
Common choices are Krylov space methods with multigrid or domain decomposition preconditioners.
Multigrid methods for adaptively refined discretizations are possible; an adaptation to isogeometric discretizations based on THB splines can be found in \cite{Hofreither2022}.
We focus on domain decomposition methods instead, where the computational domain is divided into multiple subdomains in a way that the problems local to the subdomains can be solved independently and, potentially, in parallel. These methods are especially suitable for large-scale problems and parallel computing environments. Among them, the Finite Element Tearing and Interconnecting (FETI) method, introduced in the early 1990s by Farhat and Roux in \cite{Farhat1991}, and FETI-DP, its dual-primal version, have become one of the most fundamental non-overlapping domain decomposition techniques, see \cite{Toselli2005, Mandel2005}. FETI methods enforce continuity across subdomain interfaces using Lagrange multipliers, allowing for the parallel solution of local problems while ensuring global consistency.
FETI-DP has also been adapted to multi-patch IgA in \cite{Kleiss2012} as the dual-primal Isogeometric Tearing and Interconnecting (IETI-DP) method. In IETI-DP, the patches are used as subdomains, so the patch local problem preserves the tensor-product structure that results from the use of tensor-product bases used to represent the solution on the individual patches. In~\cite{Schneckenleitner2020}, a condition number bound was given that is explicit in the discretization parameters, like the grid sizes and the spline degree. Similar results were possible for the BDDC method, see~\cite{Pavarino2013, Widlund2014}.

Those solvers were constructed for fully matching discretizations. In~\cite{Schneckenleitner2022}, the IETI-DP solver and its analysis were extended to more general geometric settings, where a discontinuous Galerkin approach \cite{Riviere2008,LangerToulopoulos2015} is used to handle the interfaces. 
In this work, we make use of the fact that the construction of adaptive splines as described in \cite{Tyoler2025} leads to spline spaces that are nested at the interfaces. We construct a IETI-DP system that does not use discontinuous Galerkin or Nitsche-type coupling. We provide a preconditioner and prove that the condition number of the preconditioned system is bounded by
$$p\left(1+\log p+\max_k \log \tfrac{H_k}{h_k}\right)^2,$$
where $p$ is the spline degree, $H_k$ is the diameter of the $k$-th patch, and $h_k$ the grid size on that patch.
That bound coincides with the bound found in~\cite{Schneckenleitner2020} for the fully matching case.

The remainder of this work is organized as follows. In Section \ref{sec:preliminaries}, we provide some basic concepts and establish the notation used throughout the paper, including the geometric setup and the multi-patch domain structure. Section \ref{sec:ieti-dp} is devoted to the derivation of the IETI-DP system and an associated preconditioner.
In Section \ref{sec:theory}, we establish the condition number bound for the preconditioned system. Finally, in Section \ref{sec:num}, we provide our numerical experiments. These experiments include an adaptive mesh refinement strategy, and we elaborate on additional challenges and remedies observed in specific problem scenarios.

\section{Model problem and preliminaries}\label{sec:preliminaries}
Let the computational domain $\Omega \subset \mathbb R^2$ be open, connected, and bounded with a Lipschitz continuous boundary $\partial \Omega$. Given a uniformly bounded diffusion parameter $\nu \in L^\infty(\Omega)$ with $\nu \geq \underline{\nu}>0$ and a right-hand-side function $f\in L^2(\Omega)$, the problem reads in variational form as follows.
\begin{equation*}\label{eq:problem}
		\mbox{Find } u\in H_0^1(\Omega):\quad
		\underbrace{\int_\Omega \nu \nabla u \cdot \nabla v \mathrm d x}_{\displaystyle a(u,v)}
		= \underbrace{\int_\Omega f v\, \mathrm d x}_{\displaystyle \langle f,v \rangle}
		\quad
		\forall v\in H^1_0(\Omega).
\end{equation*}
Here and in what follows, $H^k$, $H^k_0$ and $L^2$ are the standard Sobolev and Lebesgue spaces with standard norms $\|\cdot\|_{H^k}$ and $\|\cdot\|_{L^2}$, see~\cite{Adams2003}. To discretize this problem, we use multi-patch Isogeometric Analysis. Specifically, we assume that the computational domain $\Omega$ can be decomposed into $K$ non-overlapping subdomains (or patches) $\Omega_k$, i.e,
\begin{align*}
    \cl\Omega = \bigcup_{k=1}^K \cl\Omega_k \quad \mbox{with} \quad \Omega_k \cap \Omega_\ell = \emptyset \quad \mbox{for } k\neq\ell.
\end{align*}
Additionally, we assume that each patch $\Omega_k$ can be parameterized by a geometry function $G_k: \wh\Omega \to \Omega_k$ with parameter domain $\wh\Omega \coloneqq (0,1)^2$.

\begin{Assumption}\label{Ass:GeoEquiv}
    There exists a constant $C_G > 0$ such that
    \begin{equation*}
        \|\nabla G_k\|_{L^{\infty}(\wh\Omega)} \leq C_G H_k \quad \mbox{and} \quad \|(\nabla G_k)^{-1}\|_{L^{\infty}(\wh\Omega)} \leq C_G H^{-1}_k
    \end{equation*}
    holds for all $k = 1,\ldots,K$, where $H_k$ is the patch diameter.
\end{Assumption}

\begin{figure}[ht]
\begin{subfigure}[t]{0.45\textwidth}
    \centering
    \includegraphics[width=0.8\linewidth]{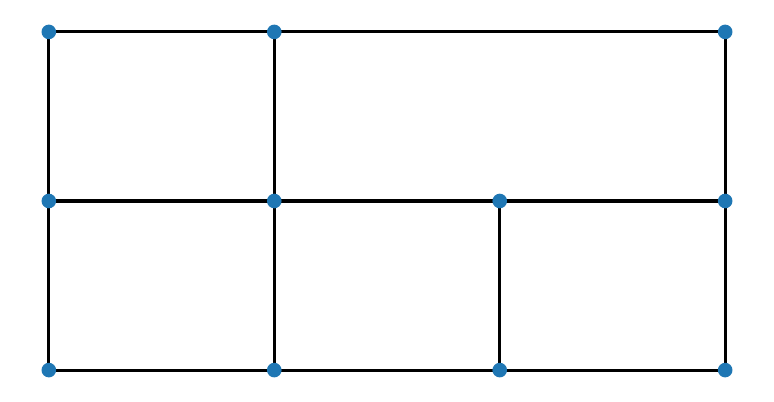}\\
    \caption{admissible}
\end{subfigure}
\hfill
\begin{subfigure}[t]{0.45\textwidth}
    \centering
    \includegraphics[width=0.8\linewidth]{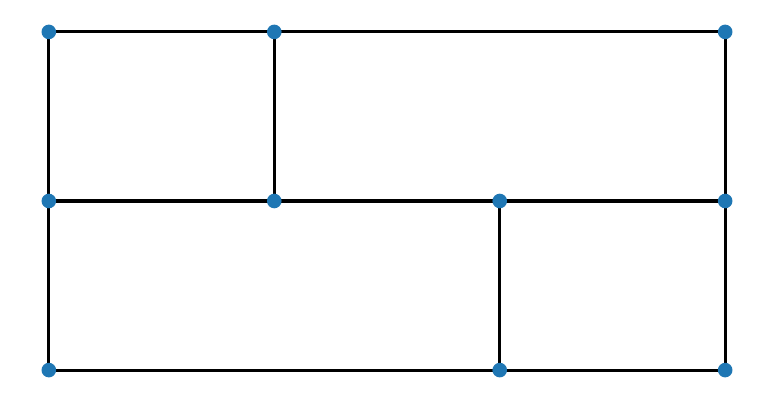}\\
    \caption{non-admissible}
\end{subfigure}
\caption{\label{fig:1}Possible multi-patch constructions}
\end{figure}

The next assumption guarantees the admissibility of the mesh in the sense of a conforming Galerkin method, see Figure~\ref{fig:1} for a visualization.
\begin{Assumption}\label{ass:geoadmissible}
    The intersection of two different patches $\cl\Omega_k \cap \cl\Omega_\ell$ is either (1) empty, (2) a common corner vertex, or (3) a whole edge of at least one of the patches.
\end{Assumption}

For more general patch decompositions, we refer to, e.g., discontinuous Galerkin (dG) methods, see \cite{Riviere2008} for an introduction to dG-methods and \cite{LangerToulopoulos2015} for its adaptation to IgA. For IETI-DP methods for such discretizations, we refer to \cite{Hofer2017,Schneckenleitner2022}. An alternative would be isogeometric mortaring methods, see~\cite{Brivadis2015}.

\medskip

For any two neighboring patches $\Omega_k$ and $\Omega_\ell$, we denote the shared edge by $\Gamma^{(k,\ell)} \coloneqq \partial\Omega_k \cap \partial\Omega_\ell$. Moreover, for every patch $\Omega_k$ 
\[
    N_\Gamma(k) \coloneqq \left\{\ell \in \{1,\ldots,K\}\,:\, \Omega_k \mbox{ and } \Omega_\ell \mbox{ share an edge}\right\}
\]
refers to its neighbors.
We denote the set of all patch-local corners and global vertices by 
\[
    \VV(\Omega_k) \coloneqq \left\{G_k(\wh\x)\,:\,\wh\x \in \{0,1\}^2\right\}\quad \mbox{and} \quad \VV \coloneqq \bigcup_{k=1}^K \VV(\Omega_k).
\]

We introduce a space of spline functions of degree $p$ for each patch $\Omega_k$. We assume to have a $p$-open knot vector in every spatial axis $\delta \in \{1,2\}$, i.e., $\Xi^{(k,\delta)} = (\xi_i^{(k,\delta)})_{i=1}^{n^{(k,\delta)}+p+1}$ such that
\begin{gather*}
    0=\xi_1^{(k,\delta)} = \cdots = \xi_{p+1}^{(k,\delta)} \le \xi_{p+2}^{(k,\delta)} \le \cdots \le \xi_{n^{(k,\delta)}}^{(k,\delta)} \le \xi_{n^{(k,\delta)}+1}^{(k,\delta)} = \cdots = \xi_{n^{(k,\delta)}+p+1}^{(k,\delta)}=1,\mbox{ and}\\
    \xi_{i}^{(k,\delta)} < \xi_{i+p}^{(k,\delta)} \quad\mbox{for}\quad i=2,\ldots,n^{(k,\delta)}.
\end{gather*}
The associated $B$-spline basis $(B^{(k,\delta)}_{i})_{i=1}^{n^{(k,\delta)}}$ is constructed by the Cox-de Boor formula, see \cite{deBoor1972}. Each basis spans a corresponding $B$-spline function space $\wh V_h^{(k,\delta)}:= \text{span}\{B^{(k,\delta)}_{i}:i=1\ldots,n^{(k,\delta)}\}\in H^1(0,1)$. For the parameter domain $\widehat \Omega$, we construct a corresponding basis $\wh\Phi^{(k)}$ with $n^{(k)} \coloneqq |\wh\Phi^{(k)}|$ as a subset of the standard tensor product B-spline basis, where we remove functions with contributions to the parametric Dirichlet boundary $\wh \Gamma_D^{(k)} \coloneqq G_k^{-1}(\partial\Omega_k \cap \partial\Omega)$. Additionally, we order the $n^{(k)}=n^{(k)}_\Ir+n^{(k)}_\Gamma$ basis functions in such a way that the first $n^{(k)}_\Ir$ basis functions only have support in the interior of $\wh\Omega$ and the remaining $n^{(k)}_\Gamma$ basis functions contribute to the boundary $\partial\wh\Omega$. Therefore we have $\wh\Phi^{(k)} = (\wh\varphi^{(k)}_i)_{i=1}^{n{(k)}}$, where
\[
    \wh\varphi_i^{(k)} \in \wh\Phi^{(k)} \iff \exists i_1,i_2 \quad \mbox{ such that } \quad \wh\varphi_i^{(k)}(\xi_1,\xi_2) = B_{i_1}^{(k,1)}(\xi_1)B_{i_2}^{(k,2)}(\xi_2) \quad \mbox{ and } \quad \wh\varphi^{(k)}_i|_{\wh\Gamma_D}=0,
\]
and furthermore
\begin{equation}\label{eq:bf order}
    \wh\varphi^{(k)}_i|_{\partial\wh\Omega}=0 \iff i \in \left\{1,\ldots,n^{(k)}_\Ir\right\} \quad \mbox{ and } \quad \wh\varphi^{(k)}_i|_{\partial\wh\Omega} \neq 0 \iff i \in n^{(k)}_\Ir + \left\{1,\ldots,n^{(k)}_\Gamma\right\}.
\end{equation}
A corresponding basis that is defined on the physical domain $\Omega_k$ is defined via the \textit{pull-back principle}, i.e.,
$\Phi^{(k)} :=(\varphi_i^{(k)})_{i=1}^{n^{(k)}}$, where $\varphi_i^{(k)} \circ\, G_k = \wh\varphi_i^{(k)}$ for $i = 1,\cdots,n^{(k)}$. These bases span the spline spaces \[
    \wh V_h^{(k)} = \mbox{span} \left\{\wh \varphi_i^{(k)} \in \wh\Phi^{(k)}\,:\,i=1,\ldots,n^{(k)}\right\} \quad \mbox{and} \quad V_h^{(k)} = \mbox{span} \left\{\varphi_i^{(k)} \in \Phi^{(k)}\,:\,i=1,\ldots,n^{(k)}\right\}.
\] 
We further define the discontinuous patch-local function space 
\[
    V_h = \prod_{k=1}^K V_h^{(k)} \quad \mbox{with} \quad n \coloneqq \text{dim}(V_h) = \sum_{k=1}^K n^{(k)}.
\]
For each patch, we define the grid size $\wh h_k$ and the smallest knot span $\wh h_{k,\min}$ via
\[
	\wh h_k := \max_{\delta\in\{1,2\}} \max_{i\in\{1,\ldots,n^{(k,\delta)}\}} \xi_{i+1}^{(k,\delta)}-\xi_{i}^{(k,\delta)},
    \quad
    \wh h_{k,\min} := \min_{\delta\in\{1,2\}} \min_{i\in\{1,\ldots,n^{(k,\delta)}\},\xi_{i+1}^{(k,\delta)}-\xi_{i}^{(k,\delta)}\ne0} \xi_{i+1}^{(k,\delta)}-\xi_{i}^{(k,\delta)}.
\]
The grid size associated with the physical domain is given by $h_k := H_k \wh h_k$.
We assume as follows.
\begin{Assumption}\label{Ass:QuasiUniform}
    There is a uniform constant $C_Q>0$ such that $\wh h_{k,\min} \ge C_Q \wh h_{k}$ holds for all $k=1,\dots,K$.
\end{Assumption}

Additionally, we define the restriction of the basis $\Phi^{(k)}$ to the edges and vertices of $\Omega_k$. Since B-splines interpolate at the boundary, the trace space on every edge is exactly spanned by the restriction of the basis functions to that edge. For all $\Gamma^{(k,\ell)} = \partial\Omega_k \cap \partial\Omega_\ell$, we define the index sets
\begin{equation*}
    \Ic^{(k)}_{\Gamma^{(k,\ell)}} \coloneqq \left\{i \in \{1,\cdots,n^{(k)}\}\,:\,\varphi_i^{(k)} \in \Phi^{(k)} \mbox{ and } \varphi_i^{(k)}|_{\Gamma^{(k,\ell)}} \not\equiv 0\right\} \quad \mbox{and}\quad \Ic^{(k)}_\Gamma \coloneqq \bigcup_{\ell \in N_\Gamma(k)} \Ic^{(k)}_{\Gamma^{(k,\ell)}}.
\end{equation*}
For all vertices $\x \in \VV \cap \cl\Omega_k$, we define the index sets
\begin{equation*}
    \Ic^{(k)}_\x \coloneqq \left\{i \in \{1,\cdots,n^{(k)}\}\,:\,\varphi_i^{(k)} \in \Phi^{(k)} \mbox{ and } \varphi_i^{(k)}(\x) \not\equiv 0\right\} \quad \mbox{and} \quad \Ic^{(k)}_{\Xc} \coloneqq \bigcup_{\x \in \VV \cap \cl\Omega_k} \Ic^{(k)}_{\x}.
\end{equation*}

The following assumption is vital in order to apply a continuous Galerkin method.
\begin{Assumption}[Nested interface spaces]\label{ass:nested}
    For any neighbouring two patches $\Omega_k$ and $\Omega_\ell$ with a shared interface $\Gamma^{(k,\ell)} = \partial \Omega_k \cap \partial\Omega_\ell$ we assume
    \[V_h^{(\ell)}|_{\Gamma^{(k,\ell)}} \subseteq V_h^{(k)}|_{\Gamma^{(k,\ell)}} \quad \mbox{or} \quad V_h^{(\ell)}|_{\Gamma^{(k,\ell)}} \supseteq V_h^{(k)}|_{\Gamma^{(k,\ell)}}.\]
\end{Assumption}
This assumption guarantees that the global discretization inherits the local approximation power of the local discretization along the interfaces between patches, see \cite{Tyoler2025} for more details. Dropping this assumption would lead to discretizations having the approximation power of the intersection of both trace spaces. Under Assumption \ref{ass:nested} we also define the following relation:
\begin{Definition}
    Let $\ell,k \in \{1,\ldots,K\}$. If $\Omega_k$ and $\Omega_\ell$ share an edge $\Gamma^{(k,\ell)} = \partial\Omega_k\cap\partial\Omega_\ell$ we define
    \[ \ell \prec k \iff V_h^{(\ell)}|_{\Gamma^{(k,\ell)}} \subset V_h^{(k)}|_{\Gamma^{(k,\ell)}} \lor \left( V_h^{(\ell)}|_{\Gamma^{(k,\ell)}} = V_h^{(k)}|_{\Gamma^{(k,\ell)}} \land \ell < k \right).\]
\end{Definition}
For any two patches sharing an edge, we have either $\ell \prec k$ or $k\prec \ell$; patches not sharing an edge are not ordered. The ordering is separate for each interface and does not imply any transitivity.
Using Assumption \ref{ass:nested}, we can further relate coarse and refined basis functions along interfaces. Again, let $\Gamma^{(k,\ell)}$ be a shared edge between $\Omega_k$ and $\Omega_\ell$ with $\ell \prec k$. Then we know that $\Phi^{(\ell)}|_{\Gamma^{(k,\ell)}} \subseteq V_h^{(k)}|_{\Gamma^{(k,\ell)}}$ and obtain
\begin{equation}\label{eq:embedding}
    \varphi_j^{(\ell)}|_{\Gamma^{(k,\ell)}} = \sum_{i \in \Ic^{(k)}_{\Gamma^{(k,\ell)}}} E^{(k,\ell)}_{i,j} \varphi_i^{(k)}|_{\Gamma^{(k,\ell)}} \quad \mbox{ for all } j \in \Ic^{(\ell)}_{\Gamma^{(k,\ell)}}.
\end{equation}
The coefficients $E^{(k,\ell)}_{i,j}$ are non-negative and are computed by a knot insertion algorithm, see \cite{deBoor1972}. Imposing $u_h^{(k)}|_{\Gamma^{(k,\ell)}} = u_h^{(\ell)}|_{\Gamma^{(k,\ell)}}$ on each interface yields the  following \textit{constraints} on the coefficients 
\begin{equation}\label{eq:ecsetup}
    u_i^{(k)} - \sum_{j \in \Ic^{(\ell)}_{\Gamma^{(k,\ell)}}}      E_{i,j}^{(k,\ell)} u_j^{(\ell)} = 0
    \quad
    \mbox{for all $\ell \prec k$ and all $i \in \Ic^{(k)}_{\Gamma^{(k,\ell)}}$.}
\end{equation}
Here and in what follows, for each discrete function, like $u_h^{(k)}$, the underlined quantity, here $\uu_h^{(k)}=(u_1^{(k)},\dots,u_{n^{(k)}}^{(k)})^\top$, is the corresponding coefficient vector.
All constraints of the form~\eqref{eq:ecsetup} are collected as rows in constraint matrices $\widetilde B_k \in \RR^{n_\lambda \times n^{(k)}}$
and
$\widetilde B = (\widetilde B_1,\ldots, \widetilde B_K) \in \RR^{n_\lambda \times N}$
such that \eqref{eq:ecsetup} is equivalent to
\begin{equation}\label{eq:nullspace}
    \widetilde B \uu_h = 0,
\end{equation}
where $n_\lambda \coloneqq \sum_{\ell \prec k} |\Ic^{(k)}_{\Gamma^{(k,\ell)}}|$ is the total number of constraints. By construction, each row of $\widetilde B=(\widetilde b_{i,j})_{1\le i\le n_\lambda, 1\le j \le n}$ contains one positive entry with value $1$ and non-positive entries otherwise, i.e, for each row index $r$ there exists a column index $j_r$ such that
\begin{equation}\label{eq:onepositive}
    \widetilde b_{r,j_r} = 1 \quad \mbox{ and } \quad \widetilde b_{r,j} \leq 0 \quad \mbox{ for all } j \neq j_r.
\end{equation}

The classical conforming Galerkin discretization reads as
\begin{equation}\label{eq:discretized}
		\mbox{Find } u_h\in V^C_h \coloneqq V_h \cap H_0^1(\Omega):\quad
		a(u_h,v_h) = \langle f,v_h\rangle
		\quad
		\forall\; v_h\in V^C_h.
\end{equation}
Using C\'ea's Lemma and the results from~\cite{Tyoler2025}, we obtain the a priori discretization error estimate
\[
    \|u-u_h\|_{H^1(\Omega)}^2
    \le c
    \inf_{v_h \in V^C_h}  \|u-v_h\|_{H^1(\Omega)}^2
    \le C
    \sum_{k=1}^K h_k^{2q_k} \|u\|_{H^{1+q_k}(\Omega_k)}^2
\]
if $u|_{\Omega_k} \in H^{1+q_k}(\Omega_k)$ and $q_k \in [1,p]$ for $k=1,\ldots,K$. A slightly weaker estimate was also proven for the case $q_k \in (0,p]$.

\section{IETI-DP}\label{sec:ieti-dp}
In this section, we derive the IETI-DP solver for the variational problem \eqref{eq:discretized} from Section \ref{sec:preliminaries}. Note that we can rewrite the bilinear form $a=\sum_{k=1}^K a^{(k)}$ and the right-hand side $f=\sum_{k=1}^K f^{(k)}$ in terms of the local contributions
\begin{gather*}
    a^{(k)}(u,v) \coloneqq \int_{\Omega_k} \nu^{(k)} \nabla u \cdot \nabla v \diff x\
    \quad\mbox{and}\quad
    \langle f^{(k)},v \rangle \coloneqq \int_{\Omega_k} f  v \diff x.
\end{gather*}
Utilizing the local bases $\Phi^{(k)}=(\varphi_i^{(k)})_{i=1}^{n^{(k)}}$, we obtain the stiffness matrix $A^{(k)} = [a^{(k)}(\varphi^{(k)}_j,\varphi^{(k)}_i)]_{i,j=1}^{n^{(k)}}$ and the local load vector $\uf^{(k)} = [\langle f,\varphi^{(k)}_i \rangle]_{i=1}^{n^{(k)}}$. Since we already established that by enforcing the continuity condition by imposing \eqref{eq:nullspace}, we see that \eqref{eq:discretized} is equivalent to the following saddle-point formulation: 
Find $(\uu_h^{(1)},\ldots,\uu_h^{(K)},\ulam)$ such that
\begin{equation}\label{eq:saddlepoint}
    \begin{pmatrix}
        A^{(1)} &         &        &         & \widetilde B^{(1)\top} \\
                & A^{(2)} &        &         & \widetilde  B^{(2)\top} \\
                &         & \ddots &         & \vdots      \\
                &         &        & A^{(K)} & \widetilde B^{(K)\top} \\
        \widetilde B^{(1)} & \widetilde B^{(2)} & \cdots & \widetilde B^{(K)} &
    \end{pmatrix}
    \begin{pmatrix}
        \uu_h^{(1)}    \\ 
        \vdots \\ 
        \uu_h^{(K)}   \\ 
        \ulam
    \end{pmatrix}
    = 
    \begin{pmatrix}
        \uf^{(1)}    \\ 
        \vdots \\ 
        \uf^{(K)}    \\ 
        0
    \end{pmatrix}.
\end{equation}

Note that the local stiffness matrices $A^{(k)}$ are singular unless $\Omega_k$ contributes to the (Dirichlet) boundary. In order to decompose the full system into sub-problems on each patch, it is crucial to resolve this issue.

\subsection{The primal space}

To guarantee the non-singularity of the local matrices, we appoint \textit{primal} degrees of freedom, i.e., we impose additional \textit{primal} constraints on the system in order to obtain a unique solution.
In the following, we split our constraints into either belonging to a vertex or not. So, we define
\[
    \Rc^{(k)}_\Xc \coloneqq \biggl\{r \in \{1,\ldots,n_\lambda\}\,:\,\exists i \in \Ic^{(k)}_{\Xc} \mbox{ such that } b^{(k)}_{r,i} = 1\biggr\} \quad \mbox{ and }\quad \Rc_\Xc \coloneqq \bigcup_{k=1}^K \Rc_\Xc^{(k)}.
\]
The matrix $\widetilde C$ contains the rows of $\widetilde B$ belonging to a vertex, $B$ the remaining ones, i.e., we define
\[B = (\widetilde b_{i,j})_{i \notin \Rc_{\Xc},1<j<n}\quad\mbox{and}\quad \widetilde C = (\widetilde b_{i,j})_{i \in \Rc_{\Xc},1<j<n}\]
The constraints in $B$ are handled using the Lagrange multiplier technique, while the constraints in $\widetilde C$ are primal. We define the vertex-continuous space
\begin{equation*}\label{eq:vertexspace}
    \widetilde V_h:=\{ v_h \in V_h : \widetilde C \uv_h = 0\}.
\end{equation*}
Additionally, define the space that satisfies the primal constraints homogeneously, i.e., splines that vanish at all vertices. Similar to the standard constraint matrix we use the splitting $\widetilde C = (\widetilde C^{(1)},\ldots, \widetilde C^{(K)})$ and introduce primal constraint matrices $C^{(k)} \in \RR^{n^{(k)}_\Pi \times n^{(k)}}$ which consist of the nonempty rows of $\widetilde C^{(k)}$. We define the global dual space as $\widetilde V_{h,\Delta} \subset \widetilde V_h$ by
\begin{equation*}\label{eq:dualspace}
    \widetilde V_{h,\Delta}  = \prod_{k=1}^K \widetilde V_{h,\Delta}^{(k)}
    \quad\mbox{with}\quad
    \widetilde V_{h,\Delta}^{(k)}:=\{ v_h \in V_h^{(k)} : C^{(k)} \uv_h^{(k)} = 0\}.
\end{equation*}
Furthermore, we define the \textit{primal space}, that is $a$-orthogonal in $\widetilde V_h$ to the dual space, i.e.,
\begin{gather*}
    \widetilde V_{h,\Pi} \coloneqq \left\{u_h \in \widetilde V_h\,:\, a^{(k)}(u_h,v_h) = 0 \mbox{ for all } k=1,\ldots,K \mbox{ and for all $v_h \in \widetilde V_{h,\Delta}$} \right\} \label{eq:primalspace}
\end{gather*}
and observe $\widetilde V_h = \widetilde V_{h,\Delta} \oplus \widetilde V_{h,\Pi}$.
We define a basis for this primal space that is nodal at all vertices. The $n_\Pi$ primal basis functions are represented by a matrix $\Psi \in \RR^{n \times n_\Pi}$, where every column gives the coefficients of the primal basis function in terms of the patch-local bases $\Phi^{(k)}$. Overall, we need to find for every patch $\Omega_k$ the basis $\Psi^{(k)} \in \RR^{n^{(k)} \times n^{(k)}_\Pi} $ such that
\begin{equation}\label{eq:localprimal}
    \uu_h^{(k)\top}A^{(k)} \Psi^{(k)} = 0 \quad \mbox{for all $\uu_h^{(k)} \in \text{ker}\,C^{(k)}$}  \quad \mbox{and} \quad C^{(k)} \Psi^{(k)} = \Ir.
\end{equation}
Since $C^{(k)}\uu_h^{(k)} = 0$, it is sufficient to have $A^{(k)}\Psi^{(k)} \in \text{Im}\,C^{(k)\top}$ and we can reformulate \eqref{eq:localprimal} in the equivalent saddle-point form
\begin{equation*}\label{eq:localprimalsaddlepoint}
    \begin{pmatrix}
        \, A^{(k)} & C^{(k)\top} \, \\
        \, C^{(k)}             \,
    \end{pmatrix}
    \begin{pmatrix}
        \, \Psi^{(k)}   \, \\
        \, \Delta^{(k)} \,
    \end{pmatrix}
    = 
    \begin{pmatrix}
        \, 0   \,\\
        \, \Ir \,
    \end{pmatrix}.
\end{equation*}
Solvability of this saddle-point system is crucial. In the absence of T-junctions, the rows of $C^{(k)}$ are linearly independent by construction since for every patch we have at most $4$ primal constraints concerning only the corner basis functions. However, for the general case, it is imperative to ensure that we do not enforce too many primal constraints (T-junctions) along any interface of the coarse patch, as the local saddle point system would be overdetermined. Thus, we require the constraint matrix (or collocation matrix) to have full row rank. In \cite{deBoor1978}, it was shown that the collocation matrix has full row rank if there is a distinct basis function $B_j$ for every collocation node $x_i$ such that $B_j(x_i)>0$. This is known as the Sch\"onberg-Whitney condition. This gives rise to the following assumption.

\begin{Assumption}[Admissibility]\label{ass:primadmissible}
    For each patch $\Omega_k$ and each of its interface edges $\Gamma^{(k,\ell)}$, the following condition holds. For every $x \in \VV \cap \cl{\Gamma^{(k,\ell)}}$ there exists a distinct $j \in \Ic^{(k)}_{\Gamma^{(k,\ell)}}$ such that $\varphi^{(k)}_j(x)>0$.
\end{Assumption}
\begin{Remark}
    In \cite{Tyoler2025}, a similar condition was needed in order to show appropriate a priori error estimates for non-matching multi-patch configurations.
\end{Remark}
The global primal basis $\Psi$ can be constructed by introducing appropriate binary restriction matrices $R^{(k)} \in \RR^{n^{(k)}_\Pi \times n_\Pi}$ that connect local primal degrees of freedom to global primal degrees of freedom 
\begin{equation*}\label{eq:globalprimal}
    \Psi
    \coloneqq 
    \begin{pmatrix}
        \; \Psi^{(1)} R^{(1)} \; \\
        \; \vspace{0.5em}\vdots \; \\
        \; \Psi^{(K)} R^{(K)} \;
    \end{pmatrix}.
\end{equation*}
By introducing the primal space, we obtain the following linear system, which is equivalent to \eqref{eq:saddlepoint}: 
Find $(\widetilde \uu_h^{(1)}, \ldots,\widetilde \uu_h^{(K)},\uu_{h,\Pi},\ulam)$ such that
\begin{equation}\label{eq:newsaddlepoint}
    \begin{pmatrix}
        \widetilde A^{(1)} &        &                    &                  & \widetilde B^{(1)\top} \\
                           & \ddots &                    &                  & \vdots                 \\
                           &        & \widetilde A^{(K)} &                  & \widetilde B^{(K)\top} \\
                           &        &                    & A_\Pi & B_{\Pi}^\top       \\
        \widetilde B^{(1)} & \cdots & \widetilde B^{(K)} & B_{\Pi}          &
    \end{pmatrix}
    \begin{pmatrix}
        \widetilde \uu_h^{(1)} \\
        \vdots \\
        \widetilde \uu_h^{(K)} \\
        \uu_{h,\Pi} \\
        \ulam
    \end{pmatrix}
    =
    \begin{pmatrix}
        \widetilde \uf^{(1)} \\
        \vdots \\
        \widetilde \uf^{(K)} \\
        \Psi^\top \uf \\
        0
    \end{pmatrix},
\end{equation}
with local system matrices $\widetilde A^{(k)} \coloneqq \begin{pmatrix}
    A^{(k)} & C^{(k)\top} \\
    C^{(k)} &              
\end{pmatrix}$ and
$A_\Pi \coloneqq \Psi^\top A \Psi$,
with local constraint matrices
$\widetilde B^{(k)} \coloneqq \begin{pmatrix}
    B^{(k)} & 0
\end{pmatrix}$ 
and
$B_{\Pi} \coloneqq B \Psi$,
with local solution vector
$\widetilde \uu_h^{(k)} \coloneqq (\uu^{(k)\top}_{h,\Delta} , \umu^{(k)\top})^\top$
and local right-hand side $\widetilde \uf^{(k)} \coloneqq (\uf^{(k)\top},0)^\top$.
The solution to the original problem is reconstructed by 
\begin{equation*}\label{eq:constructfullu}
    \uu_h^{(k)} = \uu^{(k)}_{h,\Delta} + \Psi^{(k)} \uu_{h,\Pi}.
\end{equation*}
This construction guarantees that the matrices $\widetilde A^{(k)}$ and $\Psi^\top A \Psi$   are non-singular.
By building the Schur complement for the blocks in \eqref{eq:newsaddlepoint}, we obtain the system
\begin{equation}\label{eq:schur}
    F \ulam = \ud, 
\end{equation}
where
\begin{gather}\label{eq:F}
    F \coloneqq \sum_{k=1}^K 
    \widetilde B^{(k)}
    \widetilde A^{(k)-1}
    \widetilde B^{(k)\top}
    +
    B_{\Pi} A_\Pi^{-1} B_{\Pi}^{\top}
    \quad\mbox{and}\quad 
    \ud \coloneqq \sum_{k=1}^K 
    \widetilde B^{(k)}
    \widetilde A^{(k)-1}
    \widetilde \uf^{(k)}
    + B_{\Pi} A_\Pi^{-1} \Psi^\top \uf.
\end{gather}
Observe that $F$ is symmetric positive definite.
So, in order to solve \eqref{eq:schur}, we apply a preconditioned conjugate gradient solver.

\subsection{Skeleton representation}

Following the order imposed in~\eqref{eq:bf order}, we can write the local stiffness matrices $A^{(k)}$, the local constraint matrices $B^{(k)}$ and local load vectors $\uf^{(k)}$ as
\[
    A^{(k)} = \begin{pmatrix}
        A^{(k)}_{\Ir \Ir} & A^{(k)}_{\Ir\Gamma} \\
        A^{(k)}_{\Gamma \Ir} & A^{(k)}_{\Gamma\Gamma} \\
    \end{pmatrix},
    \quad
    B^{(k)}=\begin{pmatrix} 0 & B^{(k)}_{\Gamma}\end{pmatrix}
    \quad \mbox{ and } \quad
     \uf^{(k)} = \begin{pmatrix}
        \uf^{(k)}_\Ir\\
        \uf^{(k)}_\Gamma
    \end{pmatrix}.
\]
Consequently, we can build the Schur complement with respect to interior degrees of freedom:
\begin{equation*}\label{eq:schurlocal}
    S^{(k)} \coloneqq A^{(k)}_{\Gamma\Gamma} - A^{(k)}_{\Gamma \Ir}A^{(k)-1}_{\Ir\Ir}A^{(k)}_{\Ir\Gamma} \quad \mbox{and} \quad \ug^{(k)} \coloneq \uf^{(k)}_\Gamma - A^{(k)}_{\Gamma \Ir}A^{(k)-1}_{\Ir\Ir} \uf^{(k)}_\Ir.
\end{equation*}
Using this transformation of the degrees of freedom, we can eliminate the interior degrees of freedom from any local system in terms of their skeleton degrees of freedom $\uw_h^{(k)}$.  We regain the whole solution $\uu_h^{(k)}$ by computing
\begin{equation}\label{eq:ufromw}
    \uu_h^{(k)} = \begin{pmatrix}
        A^{(k)-1}_{\Ir\Ir}(\uf^{(k)}_{\Ir}-A^{(k)}_{\Ir\Gamma} \uw_h^{(k)})\\
        \uw_h^{(k)}
    \end{pmatrix}.
\end{equation}

Analogously, we can formulate the respective spline spaces on the skeleton, i.e., we define the skeleton spaces by
\begin{equation*}\label{eq:skeleton}
    W_h^{(k)} \coloneqq \{v_h|_{\partial\Omega_k}\,:\,v_h \in V_h^{(k)}\} \quad \mbox{and} \quad W_h \coloneqq \prod_{k=1}^K W_h^{(k)}.
\end{equation*}
In the continuous setting, we can further introduce an equivalent transformation as in \eqref{eq:ufromw}. We define the \textit{discrete harmonic extension} $\Hc_h^{(k)} : W_h^{(k)} \to V_h^{(k)}$ such that
\begin{equation}\label{eq:harmonicextension}
    (\Hc_h^{(k)} w_h^{(k)})|_{\partial\Omega_k} = w_h^{(k)} \quad \mbox{and} \quad a(\Hc_h^{(k)} w_h^{(k)},v_h^{(k)})=0 \quad \mbox{for all} \quad v_h^{(k)} \in V_{h,0}^{(k)},
\end{equation}
where $V_{h,0}^{(k)} \coloneqq \{v_h \in V_h^{(k)}\,:\, v_h|_{\partial\Omega_k} = 0\}$. 

Also, the primal constraints can be represented in the skeleton formulation. Let $\widetilde C^{(k)} = \begin{pmatrix} 0  & \widetilde C_\Gamma^{(k)}\end{pmatrix}$ and $\widetilde C_\Gamma = \begin{pmatrix}\widetilde C_\Gamma^{(1)},\ldots, \widetilde C_\Gamma^{(K)}\end{pmatrix}$. We define the vertex-continuous space
\begin{equation}\label{eq:W}
    \widetilde W_h \coloneqq \{ w_h \in W_h : C_{\Gamma} \uw_h = 0\},
\end{equation}
the dual space
\begin{equation}\nonumber
    \widetilde W_{h,\Delta}^{(k)} \coloneqq \{ w_h^{(k)} \in W_h^{(k)} : C_{\Gamma}^{(k)} \uw_h^{(k)} = 0\}
    \quad \mbox{and} \quad
    \widetilde W_{h,\Delta} \coloneqq \prod_{k=1}^K W_{h,\Delta}^{(k)}
\end{equation}
and
the primal space
\begin{gather*}
    \widetilde W_{h,\Pi} \coloneqq \left\{w_h \in \widetilde W_h\,:\, a^{(k)}(\Hc_h^{(k)} w_h, \Hc_h^{(k)} q_h) = 0 \mbox{ for all } k=1,\ldots,K \mbox{ and for all $q_h \in \widetilde W_{h,\Delta}$} \right\}
\end{gather*}
and observe $\widetilde W_h = \widetilde W_{h.\Delta} \oplus \widetilde W_{h,\Pi}$.

As preconditioner, we propose the \textit{scaled Dirichlet preconditioner} $M_{sD}$, given by
\begin{equation}\label{eq:MsD}
    M_{sD} \coloneqq B_\Gamma D S D^\top B_\Gamma^\top,
\end{equation}
where the choice of the scaling matrix $D \in \RR^{n_\Gamma\times n_\Gamma}$ is described in the next section.

\subsection{Selection scaling}

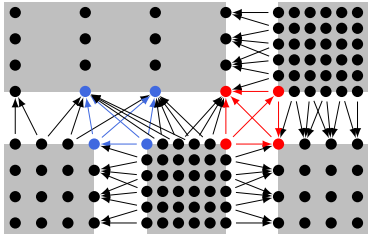
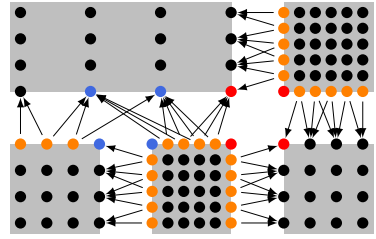
\begin{figure}[ht]
    \newcommand\dotsize{1.5}
    \newcommand\choice{orange}
    \newcommand\primone{red}
    \newcommand\primtwo{royalblue}
    \begin{subfigure}[t]{0.45\textwidth}
        \centering
        \vspace{0pt}
        \begin{tikzpicture}[scale=0.7]
	      \fill[gray!50] (-0.2,0) --   (4,0) --     (4,1.7) --    (-0.2,1.7);  
	      \fill[gray!50] (-0.2,-1) -- (1.5,-1) --   (1.5,-2.7) -- (-0.2,-2.7); 
	      \fill[gray!50] (4,-1) --    (2.5,-1) --   (2.5,-2.7) -- (4,-2.7);    
            \fill[gray!50] (5,0) --     (6.7,0) --    (6.7,1.7) --  (5,1.7);     
            \fill[gray!50] (5,-2.7) --  (6.7,-2.7) -- (6.7,-1) --   (5,-1);      

            \foreach \x in {0,0.5,1}
                \foreach \y in {-1.5,-2,-2.5}
                    \draw (\x,\y) node[circle, fill, inner sep = \dotsize pt] {};
            \draw (0,-1)     node[circle, fill, inner sep = \dotsize pt] (A1) {};
            \draw (0.5,-1)   node[circle, fill, inner sep = \dotsize pt] (A2) {};
            \draw (1,-1)     node[circle, fill, inner sep = \dotsize pt] (A3) {};
            \draw (1.5,-1)   node[circle, fill, inner sep = \dotsize pt, color=\primtwo] (A4) {};
            \draw (1.5,-1.5) node[circle, fill, inner sep = \dotsize pt] (A5) {};
            \draw (1.5,-2)   node[circle, fill, inner sep = \dotsize pt] (A6) {};
            \draw (1.5,-2.5) node[circle, fill, inner sep = \dotsize pt] (A7) {};

            \foreach \x in {2.8,3.1,3.4,3.7}
                \foreach \y in {-1.3,-1.6,-1.9,-2.2,-2.5}
                    \draw (\x,\y) node[circle, fill, inner sep = \dotsize pt] {};
            \draw (2.5,-2.5) node[circle, fill, inner sep = \dotsize pt] (B1) {};
            \draw (2.5,-2.2) node[circle, fill, inner sep = \dotsize pt] (B2) {};
            \draw (2.5,-1.9) node[circle, fill, inner sep = \dotsize pt] (B3) {};
            \draw (2.5,-1.6) node[circle, fill, inner sep = \dotsize pt] (B4) {};
            \draw (2.5,-1.3) node[circle, fill, inner sep = \dotsize pt] (B5) {};
            \draw (2.5,-1)   node[circle, fill, inner sep = \dotsize pt, color=\primtwo] (B6) {};
            \draw (2.8,-1)   node[circle, fill, inner sep = \dotsize pt] (B7) {};
            \draw (3.1,-1)   node[circle, fill, inner sep = \dotsize pt] (B8) {};
            \draw (3.4,-1)   node[circle, fill, inner sep = \dotsize pt] (B9) {};
            \draw (3.7,-1)   node[circle, fill, inner sep = \dotsize pt] (B10) {};
            \draw (4,-1)     node[circle, fill, inner sep = \dotsize pt, color=\primone] (B11) {};
            \draw (4,-1.3)   node[circle, fill, inner sep = \dotsize pt] (B12) {};
            \draw (4,-1.6)   node[circle, fill, inner sep = \dotsize pt] (B13) {};
            \draw (4,-1.9)   node[circle, fill, inner sep = \dotsize pt] (B14) {};
            \draw (4,-2.2)   node[circle, fill, inner sep = \dotsize pt] (B15) {};
            \draw (4,-2.5)   node[circle, fill, inner sep = \dotsize pt] (B16) {};

            \foreach \x in {5.5,6,6.5}
                \foreach \y in {-1.5,-2,-2.5}
                    \draw (\x,\y) node[circle, fill, inner sep = \dotsize pt] {};
            \draw (5,-2.5) node[circle, fill, inner sep = \dotsize pt] (C1) {};
            \draw (5,-2)   node[circle, fill, inner sep = \dotsize pt] (C2) {};
            \draw (5,-1.5) node[circle, fill, inner sep = \dotsize pt] (C3) {};
            \draw (5,-1)   node[circle, fill, inner sep = \dotsize pt,color=\primone] (C4) {};
            \draw (5.5,-1) node[circle, fill, inner sep = \dotsize pt] (C5) {};
            \draw (6,-1)   node[circle, fill, inner sep = \dotsize pt] (C6) {};
            \draw (6.5,-1) node[circle, fill, inner sep = \dotsize pt] (C7) {};
            
            \foreach \x in {5.3,5.6,5.9,6.2,6.5}
                \foreach \y in {1.5,1.2,0.9,0.6,0.3}
                    \draw (\x,\y) node[circle, fill, inner sep = \dotsize pt] {};
            \draw (6.5,0) node[circle, fill, inner sep = \dotsize pt] (D1) {};
            \draw (6.2,0) node[circle, fill, inner sep = \dotsize pt] (D2) {};
            \draw (5.9,0) node[circle, fill, inner sep = \dotsize pt] (D3) {};
            \draw (5.6,0) node[circle, fill, inner sep = \dotsize pt] (D4) {};
            \draw (5.3,0) node[circle, fill, inner sep = \dotsize pt] (D5) {};
            \draw (5,0)   node[circle, fill, inner sep = \dotsize pt,color=\primone] (D6) {};
            \draw (5,0.3) node[circle, fill, inner sep = \dotsize pt] (D7) {};
            \draw (5,0.6) node[circle, fill, inner sep = \dotsize pt] (D8) {};
            \draw (5,0.9) node[circle, fill, inner sep = \dotsize pt] (D9) {};
            \draw (5,1.2) node[circle, fill, inner sep = \dotsize pt] (D10) {};
            \draw (5,1.5) node[circle, fill, inner sep = \dotsize pt] (D11) {};

            \foreach \x in {0,1.33,2.66}
                \foreach \y in {0.5,1,1.5}
                    \draw (\x,\y) node[circle, fill, inner sep = \dotsize pt] {};
            \draw (4,1.5)  node[circle, fill, inner sep = \dotsize pt] (E1) {};
            \draw (4,1)    node[circle, fill, inner sep = \dotsize pt] (E2) {};
            \draw (4,0.5)  node[circle, fill, inner sep = \dotsize pt] (E3) {};
            \draw (4,0)    node[circle, fill, inner sep = \dotsize pt,color=\primone] (E4) {};
            \draw (2.66,0) node[circle, fill, inner sep = \dotsize pt, color=\primtwo] (E5) {};
            \draw (1.33,0) node[circle, fill, inner sep = \dotsize pt, color=\primtwo] (E6) {};
            \draw (0,0)    node[circle, fill, inner sep = \dotsize pt] (E7) {};

            \draw[-latex, line width = 0.3pt, shorten <=\dotsize]
            (B6) edge[color=\primtwo] (A4) (B5) edge (A4) (B5) edge (A5) (B4) edge (A5) (B4) edge (A6) (B3) edge (A5) (B3) edge (A6) (B2) edge (A6) (B2) edge (A7) (B1) edge (A7)
            (A1) edge (E7) (A2) edge (E7) (A2) edge (E6) (A3) edge (E6) (A3) edge (E5) (A4) edge[color=\primtwo] (E6) (A4) edge[color=\primtwo] (E5)
            (B16) edge (C1) (B15) edge (C1) (B15) edge (C2) (B14) edge (C2) (B14) edge (C3) (B13) edge (C2) (B13) edge (C3) (B12) edge (C3) (B12) edge (C4) (B11) edge[color=\primone] (C4) (B11) edge[color=\primone] (D6)
            (B6) edge[color=\primtwo] (E6) (B6) edge[color=\primtwo] (E5) (B7) edge (E6) (B7) edge (E5) (B8) edge (E6) (B8) edge (E5) (B8) edge (E4) (B9) edge (E6) (B9) edge (E5) (B9) edge (E4) (B10) edge (E5) (B10) edge (E4) (B11) edge[color=\primone] (E4)
            (D1) edge (C7) (D2) edge (C7) (D2) edge (C6) (D3) edge (C6) (D3) edge (C5) (D4) edge (C6) (D4) edge (C5) (D5) edge (C5) (D5) edge (C4) (D6) edge[color=\primone] (C4) (C4) edge[color=\primone] (E4)
            (D6) edge[color=\primone] (E4) (D7) edge (E4) (D7) edge (E3) (D8) edge (E3) (D8) edge (E2) (D9) edge (E3) (D9) edge (E2) (D10) edge (E2) (D10) edge (E1) (D11) -- (E1);
        \end{tikzpicture}
        \caption{Continuity constraints for non-matching tensor-product B-splines. Each arrow represents a negative entry and each degree of freedom from which an arrow emerges represents a $1$ in the constraint matrix $\widetilde B$. Colored degrees of freedom belong to $\Ic_\x^{(k)}$ and colored constraints belong to $\Rc^{(k)}_\Xc$ for different vertices $\x \in \VV$ and different patches $\Omega_k$.}
    \end{subfigure}
    \hfill
    \begin{subfigure}[t]{0.45\textwidth}
        \centering
        \vspace{0pt}
        \begin{tikzpicture}[scale=0.7]
	      \fill[gray!50] (-0.2,0) --   (4,0) --     (4,1.7) --    (-0.2,1.7);  
	      \fill[gray!50] (-0.2,-1) -- (1.5,-1) --   (1.5,-2.7) -- (-0.2,-2.7); 
	      \fill[gray!50] (4,-1) --    (2.5,-1) --   (2.5,-2.7) -- (4,-2.7);    
            \fill[gray!50] (5,0) --     (6.7,0) --    (6.7,1.7) --  (5,1.7);     
            \fill[gray!50] (5,-2.7) --  (6.7,-2.7) -- (6.7,-1) --   (5,-1);      

            \foreach \x in {0,0.5,1}
                \foreach \y in {-1.5,-2,-2.5}
                    \draw (\x,\y) node[circle, fill, inner sep = \dotsize pt] {};
            \draw (0,-1)     node[circle, fill, inner sep = \dotsize pt, color=\choice] (A1) {};
            \draw (0.5,-1)   node[circle, fill, inner sep = \dotsize pt, color=\choice] (A2) {};
            \draw (1,-1)     node[circle, fill, inner sep = \dotsize pt, color=\choice] (A3) {};
            \draw (1.5,-1)   node[circle, fill, inner sep = \dotsize pt, color=\primtwo] (A4) {};
            \draw (1.5,-1.5) node[circle, fill, inner sep = \dotsize pt] (A5) {};
            \draw (1.5,-2)   node[circle, fill, inner sep = \dotsize pt] (A6) {};
            \draw (1.5,-2.5) node[circle, fill, inner sep = \dotsize pt] (A7) {};

            \foreach \x in {2.8,3.1,3.4,3.7}
                \foreach \y in {-1.3,-1.6,-1.9,-2.2,-2.5}
                    \draw (\x,\y) node[circle, fill, inner sep = \dotsize pt] {};
            \draw (2.5,-2.5) node[circle, fill, inner sep = \dotsize pt, color=\choice] (B1) {};
            \draw (2.5,-2.2) node[circle, fill, inner sep = \dotsize pt, color=\choice] (B2) {};
            \draw (2.5,-1.9) node[circle, fill, inner sep = \dotsize pt, color=\choice] (B3) {};
            \draw (2.5,-1.6) node[circle, fill, inner sep = \dotsize pt, color=\choice] (B4) {};
            \draw (2.5,-1.3) node[circle, fill, inner sep = \dotsize pt, color=\choice] (B5) {};
            \draw (2.5,-1)   node[circle, fill, inner sep = \dotsize pt, color=\primtwo] (B6) {};
            \draw (2.8,-1)   node[circle, fill, inner sep = \dotsize pt, color=\choice] (B7) {};
            \draw (3.1,-1)   node[circle, fill, inner sep = \dotsize pt, color=\choice] (B8) {};
            \draw (3.4,-1)   node[circle, fill, inner sep = \dotsize pt, color=\choice] (B9) {};
            \draw (3.7,-1)   node[circle, fill, inner sep = \dotsize pt, color=\choice] (B10) {};
            \draw (4,-1)     node[circle, fill, inner sep = \dotsize pt, color=\primone] (B11) {};
            \draw (4,-1.3)   node[circle, fill, inner sep = \dotsize pt, color=\choice] (B12) {};
            \draw (4,-1.6)   node[circle, fill, inner sep = \dotsize pt, color=\choice] (B13) {};
            \draw (4,-1.9)   node[circle, fill, inner sep = \dotsize pt, color=\choice] (B14) {};
            \draw (4,-2.2)   node[circle, fill, inner sep = \dotsize pt, color=\choice] (B15) {};
            \draw (4,-2.5)   node[circle, fill, inner sep = \dotsize pt, color=\choice] (B16) {};

            \foreach \x in {5.5,6,6.5}
                \foreach \y in {-1.5,-2,-2.5}
                    \draw (\x,\y) node[circle, fill, inner sep = \dotsize pt] {};
            \draw (5,-2.5) node[circle, fill, inner sep = \dotsize pt] (C1) {};
            \draw (5,-2)   node[circle, fill, inner sep = \dotsize pt] (C2) {};
            \draw (5,-1.5) node[circle, fill, inner sep = \dotsize pt] (C3) {};
            \draw (5,-1)   node[circle, fill, inner sep = \dotsize pt, color=\primone] (C4) {};
            \draw (5.5,-1) node[circle, fill, inner sep = \dotsize pt] (C5) {};
            \draw (6,-1)   node[circle, fill, inner sep = \dotsize pt] (C6) {};
            \draw (6.5,-1) node[circle, fill, inner sep = \dotsize pt] (C7) {};
            
            \foreach \x in {5.3,5.6,5.9,6.2,6.5}
                \foreach \y in {1.5,1.2,0.9,0.6,0.3}
                    \draw (\x,\y) node[circle, fill, inner sep = \dotsize pt] {};
            \draw (6.5,0) node[circle, fill, inner sep = \dotsize pt, color=\choice] (D1) {};
            \draw (6.2,0) node[circle, fill, inner sep = \dotsize pt, color=\choice] (D2) {};
            \draw (5.9,0) node[circle, fill, inner sep = \dotsize pt, color=\choice] (D3) {};
            \draw (5.6,0) node[circle, fill, inner sep = \dotsize pt, color=\choice] (D4) {};
            \draw (5.3,0) node[circle, fill, inner sep = \dotsize pt, color=\choice] (D5) {};
            \draw (5,0)   node[circle, fill, inner sep = \dotsize pt, color=\primone] (D6) {};
            \draw (5,0.3) node[circle, fill, inner sep = \dotsize pt, color=\choice] (D7) {};
            \draw (5,0.6) node[circle, fill, inner sep = \dotsize pt, color=\choice] (D8) {};
            \draw (5,0.9) node[circle, fill, inner sep = \dotsize pt, color=\choice] (D9) {};
            \draw (5,1.2) node[circle, fill, inner sep = \dotsize pt, color=\choice] (D10) {};
            \draw (5,1.5) node[circle, fill, inner sep = \dotsize pt, color=\choice] (D11) {};

            \foreach \x in {0,1.33,2.66}
                \foreach \y in {0.5,1,1.5}
                    \draw (\x,\y) node[circle, fill, inner sep = \dotsize pt] {};
            \draw (4,1.5)  node[circle, fill, inner sep = \dotsize pt] (E1) {};
            \draw (4,1)    node[circle, fill, inner sep = \dotsize pt] (E2) {};
            \draw (4,0.5)  node[circle, fill, inner sep = \dotsize pt] (E3) {};
            \draw (4,0)    node[circle, fill, inner sep = \dotsize pt, color=\primone] (E4) {};
            \draw (2.66,0) node[circle, fill, inner sep = \dotsize pt, color=\primtwo] (E5) {};
            \draw (1.33,0) node[circle, fill, inner sep = \dotsize pt, color=\primtwo] (E6) {};
            \draw (0,0)    node[circle, fill, inner sep = \dotsize pt] (E7) {};

            \draw[-latex, line width = 0.3pt, shorten <=\dotsize]
            (B5) edge (A4) (B5) edge (A5) (B4) edge (A5) (B4) edge (A6) (B3) edge (A5) (B3) edge (A6) (B2) edge (A6) (B2) edge (A7) (B1) edge (A7)
            (A1) edge (E7) (A2) edge (E7) (A2) edge (E6) (A3) edge (E6) (A3) edge (E5)
            (B16) edge (C1) (B15) edge (C1) (B15) edge (C2) (B14) edge (C2) (B14) edge (C3) (B13) edge (C2) (B13) edge (C3) (B12) edge (C3) (B12) edge (C4)  
            (B7) edge (E6) (B7) edge (E5) (B8) edge (E6) (B8) edge (E5) (B8) edge (E4) (B9) edge (E6) (B9) edge (E5) (B9) edge (E4) (B10) edge (E5) (B10) edge (E4) 
            (D1) edge (C7) (D2) edge (C7) (D2) edge (C6) (D3) edge (C6) (D3) edge (C5) (D4) edge (C6) (D4) edge (C5) (D5) edge (C5) (D5) edge (C4) 
            (D7) edge (E4) (D7) edge (E3) (D8) edge (E3) (D8) edge (E2) (D9) edge (E3) (D9) edge (E2) (D10) edge (E2) (D10) edge (E1) (D11) -- (E1);
	\end{tikzpicture}
        \caption{One may remove redundant constraints from $\widetilde B$ since we already enforce continuity on all vertices in the space $\widetilde W_h$, i.e. the corresponding degrees of freedom are already joined in the system due to the primal space (by computing $\Psi$). The local degrees of freedom in $\Ic_\Zc$ are depicted in \choice. By removing corner constraints one will always obtain $\Ic_\Zc \cap \Rc_\Xc = \emptyset$.}
        \label{fig:selectionscaling}
    \end{subfigure}
    \caption{Schematic representation of constraints on non-matching interfaces}
    \label{fig:nonmatching}
\end{figure}

In order to find a suitable scaling matrix $D$, we collect constraints and local degrees of freedom that are connected via the constraint matrix $B_\Gamma$ and define the set of constraints belonging to a certain degree of freedom, i.e, for a local degree of freedom $i$ of patch $\Omega_k$ we define
\begin{equation*}
    \Rc_i^{(k)} \coloneqq \left\{r \in \{1,\dots,n_\lambda\}\,:\,b_{r,i}^{(k)} \neq 0 \right\}.
\end{equation*}
Moreover, we define the local degrees of freedom that are on the refined side of an interface as seen in Figure \ref{fig:selectionscaling}, i.e.,
\begin{equation*}
    \Ic_\Zc \coloneqq \left\{(k,j)\,:\,\exists! \ell \text{ with } \ell \prec k \mbox{ such that } j \in \Ic_{\Gamma^{(k,\ell)}}^{(k)}\right\}.
\end{equation*}
Equivalently, the degrees of freedom in $\Ic_{\Zc}$ can be characterized by the fact that there is exactly $1$ constraint (row) in the constraint matrix with the respective $j$-th entry being $1$. We define the mapping $\Zc \,:\, \Ic_\Zc \to \{1,\cdots,n_\lambda\}$ so that $\Zc$ returns exactly the specific row index corresponding to the degree of freedom $j$ of patch $\Omega_k$, where $(k,j) \in \Ic_\Zc$, that is, $\Zc : (k,j) \mapsto r \mbox{ such that } b^{(k)}_{r,j} = 1$ and $b^{(k)}_{i,j} = 0$ for $i \neq r$.
\begin{Lemma}\label{lem:bijection}
    The mapping $\Zc : \Ic_\Zc \to \{1,\dots,n_\lambda\}$ is injective.
\end{Lemma}
\begin{Proof}
    The statement follows directly from \eqref{eq:onepositive}.
\end{Proof}
\begin{Definition}[Selection scaling]\label{def:selection}
We define $D \coloneqq \text{diag}\,(D^{(1)}, \ldots, D^{(K)} ) \in \RR^{n_\Gamma \times n_\Gamma}$ as a diagonal matrix with entries being $1$ whenever the index on the diagonal refers to a degree of freedom on the refined side of an interface that corresponds to exactly one constraint. Thus, the entries $d_{i,j}^{(k)}$ of $D^{(k)} \in \RR^{n_\Gamma^{(k)}\times n_\Gamma^{(k)}}$ are given by
\begin{equation}\label{eq:D def}
    d_{i,j}^{(k)} = \begin{cases}
        1 \quad \mbox{if } i=j \mbox{ and } (k,i) \in \Ic_\Zc, \\
        0 \quad \mbox{else}.
    \end{cases}
\end{equation}
\end{Definition}

\section{Condition number estimate}\label{sec:theory}
In this section, we present the main contribution of this paper. Specifically, we show a condition number estimate for the preconditioned Schur complement system. Throughout this section, we use the following notation to indicate upper bounds without specifying certain constants:
\begin{Definition}
    We write $A \lesssim B$ if there is a generic constant $C$ that 
    only depends
    on the constant $C_G$ from Assumption~\ref{Ass:GeoEquiv}
    and
    the constant $C_Q$ from Assumption~\ref{Ass:QuasiUniform}
    such that $A\leq CB$.
\end{Definition}

Next, we refer to \cite{Mandel2005} by Mandel, Dohrmann, and Tezaur, where the following abstract error bound was proven:
\begin{Theorem}\label{thm:abstractcondition}
    Let $\widetilde W_h$ be given as in \eqref{eq:W} and let $M_{sD}$ and $F$ be given as in \eqref{eq:F} and \eqref{eq:MsD}. If for every $w_h \in \widetilde W_h$ and $v_h \in W_h$ such that $\uv_h = D B_\Gamma^\top B_\Gamma \uw_h$ it holds that 
    \begin{itemize}
        \item $v_h \in \widetilde W_h$,
        \item $B_\Gamma \uv_{h} = B_\Gamma \uw_{h}$ if additionally $w_h \in \widetilde W_{h,\Delta}$,
    \end{itemize} then we obtain the abstract condition number bound
    \begin{equation}\label{eq:omega}
        \kappa(M_{sD}F) = \frac{\lambda_{\text{max}}(M_{sD}F)}{\lambda_{\text{min}}(M_{sD}F)} \leq \sup_{0\neq w_h\in\widetilde W_h}\frac{\|D B_\Gamma^\top B_\Gamma \uw_h\|^2_S}{\|\uw_h\|^2_S}
    \end{equation}
\end{Theorem}
\begin{Proof}
   See \cite[Theorem 22]{Mandel2005}.
\end{Proof}

To apply Theorem \ref{thm:abstractcondition} we must check the two conditions for our setting of constraints $B_\Gamma$ between non-matching patches and for our selection scaling matrix $D$. We first show the following auxiliary lemma concerning the structure of the skeleton function after applying the operator $DB_\Gamma^\top B_\Gamma$.

\begin{Lemma}\label{lem:w}
    Let $w_h = (w_h^{(1)}, \cdots, w_h^{(K)}) \in \widetilde W_h$ and $v_h = (v_h^{(1)}, \cdots, v_h^{(K)}) \in W_h$ such that 
    $
        \uv_h = DB_\Gamma^\top B_\Gamma \uw_h.
    $
     Then, for each patch $\Omega_k$ and each interface $\Gamma^{(k,\ell)}$, the following identity holds:
    \begin{equation*}\label{eq:normwu}
        v_h^{(k)}|_{\Gamma^{(k,\ell)}} = \begin{cases} w_h^{(k)}|_{\Gamma^{(k,\ell)}} - w_h^{(\ell)}|_{\Gamma^{(k,\ell)}}\, &\mbox{if $\ell \prec k$}, \\ 0 \, &\mbox{else}.\end{cases}
    \end{equation*}
\end{Lemma}
\begin{Proof}
     Applying the matrix $B_\Gamma$ to $\uw_{h}$ yields the Lagrange multipliers, i.e.,
    \begin{align*}
        (B_\Gamma \uw_{h})_r = \sum_{k=1}^K \sum_{j \in \Ic_\Gamma^{(k)}} b_{r,j}^{(k)} w_{j}^{(k)} \eqqcolon \lambda_r.
    \end{align*}
    Note that if $r \in \Rc_\Xc$, we obtain $(B_\Gamma \uw_{h})_r = 0$ since $w_{h} \in \widetilde W_{h}$ and therefore the corner constraints are satisfied. Observe that the entry $v_i^{(k)}$ related to degree of freedom $i$ of patch $\Omega_k$ can be expressed by
    \[
        v_i^{(k)} = d_{i,i}^{(k)} \sum_{r \in \Rc_i^{(k)}} b_{r,i}^{(k)} \lambda_r,
    \]
    since $D$ is a diagonal matrix. Note that $d_{i,i}^{(k)}\neq 0$ if $(k,i) \in \Ic_\Zc$, i.e., the set $\Rc_i^{(k)}$ has exactly one element, namely $\Rc_i^{(k)} = \{\Zc((k,i))\}$ and additionally $b_{\Zc((k,i)),i}^{(k)} = 1$. Otherwise $d_{i,i}^{(k)} = 0$, which means we obtain the simplified expression
    \begin{align*}   
        v_i^{(k)} = \begin{cases}
           \sum_{p=1}^K \sum_{j \in \Ic_\Gamma^{(p)}} b_{\Zc((p,i)),j}^{(p)} w_j^{(p)} \quad &\mbox{if } (k,i) \in \Ic_\Zc,\\
            0 \quad &\mbox{else}.
        \end{cases}
    \end{align*}
    Let $\Gamma^{(k,\ell)} = \partial \Omega_k \cap \partial \Omega_\ell$ with $\ell \prec k$. Furthermore, we denote the endpoints of the interface by $\x^{(k,\ell,1)}$ and $\x^{(k,\ell,2)}$ respectively, and $\varphi_{i_\alpha}^{(k)}$ denotes the basis function in $\Phi^{(k)}$ such that $\varphi_{i_\alpha}^{(k)}(\x^{(k,\ell,\alpha)})=1$ for $\alpha = 1,2$. Without loss of generality, we assume that $(k,i_\alpha) \notin \Ic_\Zc$, that is, $v^{(k)}_{i_\alpha} = 0$. The corresponding spline function in patch $\Omega_k$ restricted to an interface $\Gamma^{(k,\ell)}$ can be expressed by
    \begin{align*}
        v_h^{(k)}|_{\Gamma^{(k,\ell)}} &= \sum_{i \in \Ic^{(k)}_{\Gamma^{(k,\ell)}} } \varphi_i^{(k)}|_{\Gamma^{(k,\ell)}}  v_i^{(k)} = \sum_{i\in \Ic^{(k)}_{\Gamma^{(k,\ell)}}\setminus \{i_1,i_2\} } \varphi_i^{(k)}|_{\Gamma^{(k,\ell)}}  \biggl(\sum_{p=1}^K \sum_{j \in \Ic_{\Gamma^{(k,\ell)}}^{(p)}} b_{\Zc((k,i)),j}^{(p)} w_j^{(p)}\biggr).
    \end{align*}
    Note that $(k,i) \in \Ic_\Zc$ if and only if the corresponding basis function has only one constraint, i.e, at most only $2$ patches are involved. Since we restricted ourselves to $\Gamma^{(k,\ell)}$, the only possible pair is $\Omega_k$ and $\Omega_\ell$, thus we obtain
    \begin{align*}
        v_h^{(k)}|_{\Gamma^{(k,\ell)}} = \sum_{i\in \Ic^{(k)}_{\Gamma^{(k,\ell)}}\setminus\{i_1,i_2\} } \varphi_i^{(k)}|_{\Gamma^{(k,\ell)}}  \biggl(\sum_{j \in \Ic_{\Gamma^{(k,\ell)}}^{(k)}} b_{\Zc((k,i)),j}^{(k)} w_j^{(k)} + \sum_{j \in \Ic_{\Gamma^{(k,\ell)}}^{(\ell)}} b_{\Zc((k,i)),j}^{(\ell)} w_j^{(\ell)}\biggr).
    \end{align*}
    Since $\Rc_i^{(k)} = \{\Zc((k,i))\}$ and $b_{\Zc((k,i)),i}^{(k)} = 1$, $b_{\Zc((k,i)),j}^{(k)} = 0$ must hold for $j \neq i$, so
    \begin{align*}
        v_h^{(k)}|_{\Gamma^{(k,\ell)}} &= \sum_{i\in \Ic^{(k)}_{\Gamma^{(k,\ell)}} \setminus \{i_1,i_2\}} \varphi_i^{(k)}|_{\Gamma^{(k,\ell)}}  \biggl(w_i^{(k)} + \sum_{j \in \Ic_{\Gamma^{(k,\ell)}}^{(\ell)}} b_{\Zc((k,i)),j}^{(\ell)} w_j^{(\ell)}\biggr) \\
        &= \sum_{i\in \Ic^{(k)}_{\Gamma^{(k,\ell)}}\setminus \{i_1,i_2\}} \varphi_i^{(k)}|_{\Gamma^{(k,\ell)}}  w_i^{(k)} + \sum_{j \in \Ic_{\Gamma^{(k,\ell)}}^{(\ell)}} \sum_{i\in \Ic^{(k)}_{\Gamma^{(k,\ell)}}\setminus \{i_1,i_2\}} \varphi_i^{(k)}|_{\Gamma^{(k,\ell)}}  b_{\Zc((k,i)),j}^{(\ell)}w_j^{(\ell)}.
    \end{align*}
    Observe that we obtain for the first term
    \begin{align}\label{eq:firstterm}
        \sum_{i\in \Ic^{(k)}_{\Gamma^{(k,\ell)}}\setminus \{i_1,i_2\}} \varphi_i^{(k)}|_{\Gamma^{(k,\ell)}}  w_i^{(k)} &= \underbrace{\sum_{i\in \Ic^{(k)}_{\Gamma^{(k,\ell)}}} \varphi_i^{(k)}|_{\Gamma^{(k,\ell)}}  w_i^{(k)}}_{= w^{(k)}|_{\Gamma^{k,\ell)}}} - \sum_{\alpha=1}^2 \varphi_{i_\alpha}^{(k)}|_{\Gamma^{(k,\ell)}}  w_h^{(k)}(\x^{(k,\ell,\alpha)})
    \end{align}
    and for the second term
    \begin{multline*}
        \sum_{j \in \Ic_{\Gamma^{(k,\ell)}}^{(\ell)}} \sum_{i\in \Ic^{(k)}_{\Gamma^{(k,\ell)}}\setminus \{i_1,i_2\}} \varphi_i^{(k)}|_{\Gamma^{(k,\ell)}} b_{\Zc((k,i)),j}^{(\ell)}w_j^{(\ell)} \\
        = \sum_{j \in \Ic_{\Gamma^{(k,\ell)}}^{(\ell)}} \biggl(\sum_{i\in \Ic^{(k)}_{\Gamma^{(k,\ell)}}} \varphi_i^{(k)}|_{\Gamma^{(k,\ell)}} b_{r^{(k,\ell)}_i,j}^{(\ell)} \biggr)  w_j^{(\ell)}
        - \sum_{\alpha=1}^2\varphi_{i_\alpha}^{(k)}|_{\Gamma^{(k,\ell)}} \sum_{j \in \Ic_{\Gamma^{(k,\ell)}}^{(\ell)}}b^{(\ell)}_{r^{(k,\ell)}_{i_\alpha},j} w^{(\ell)}_j ,
    \end{multline*}
    where we denote $r^{(k,\ell)}_i = \Zc((k,i))$ if $(k,i) \in \Ic_\Zc$. For the indices at the endpoints of the interface ($i_1$ and $i_2$) $\Zc$ is not well defined and therefore, we define $r^{(k,\ell)}_{i_\alpha}$ as the constraint connecting $i_\alpha$ to a related degree of freedom on patch $\Omega_\ell$ (there are in general multiple constraints to multiple patches since $i_\alpha$ is a corner degree of freedom). Furthermore, we know for the first expression that the entries of $B$ correspond to the negative embedding coefficients, and by using \eqref{eq:embedding}, we obtain
    \[\sum_{i\in \Ic^{(k)}_{\Gamma^{(k,\ell)}}} \varphi_i^{(k)}|_{\Gamma^{(k,\ell)}} b_{r^{(k,\ell)}_i,j}^{(\ell)} = -\sum_{i\in \Ic^{(k)}_{\Gamma^{(k,\ell)}}} E_{i,j}^{(k,\ell)} \varphi_i^{(k)}|_{\Gamma^{(k,\ell)}} = -\varphi_j^{(\ell)}|_{\Gamma^{(k,\ell)}} \quad \mbox{for all } j \in \Ic^{(\ell)}_{\Gamma^{(k,\ell)}}.\]
    For the corner constraints, it holds that
    \[\sum_{j \in \Ic_{\Gamma^{(k,\ell)}}^{(\ell)}}b^{(\ell)}_{r^{(k,\ell)}_{i_\alpha},j} w_j^{(\ell)} = -\sum_{j \in \Ic_{\Gamma^{(k,\ell)}}^{(\ell)}} E^{(k,\ell)}_{i_\alpha,j} w_j^{(\ell)} = -w_h^{(\ell)}(\x^{(k,\ell,\alpha)})\quad \mbox{for } \alpha=1,2,\]
    where we used the fact that the constraint matrix (or negative embedding matrix) maps the coarse coefficients exactly to the corresponding coefficient of a nodal function equal to $\varphi^{(k)}_{i_\alpha}|_{\Gamma^{(k,\ell)}}$. Since this function is nodal at $\x^{(k,\ell,\alpha)}$, the coefficient is equal to $w_h^{(\ell)}(\x^{(k,\ell,\alpha)})$. Thus, we have
    \begin{multline}\label{eq:secterm}
        \sum_{j \in \Ic_{\Gamma^{(k,\ell)}}^{(\ell)}} \sum_{i\in \Ic^{(k)}_{\Gamma^{(k,\ell)}}\setminus \{i_1,i_2\}} \varphi_i^{(k)}|_{\Gamma^{(k,\ell)}} b_{\Zc((k,i)),j}^{(\ell)}w_j^{(\ell)}\\
        = - \underbrace{\sum_{j \in \Ic_{\Gamma^{(k,\ell)}}^{(\ell)}} \varphi_j^{(j)}|_{\Gamma^{(k,\ell)}}  w_j^{(\ell)}}_{w^{(\ell)}|_{\Gamma^{(k,\ell)}}}
        + \sum_{\alpha=1}^2\varphi_{i_\alpha}^{(k)}|_{\Gamma^{(k,\ell)}} w_h^{(\ell)}(\x^{(k,\ell,\alpha)}) ,
    \end{multline}
    By combining \eqref{eq:firstterm} and \eqref{eq:secterm} we finally obtain
    \begin{align*}
        v_h^{(k)}|_{\Gamma^{(k,\ell)}} = w_h^{(k)}|_{\Gamma^{(k,\ell)}} - w_h^{(\ell)}|_{\Gamma^{(k,\ell)}} &- \sum_{\alpha=1}^2\varphi^{(k)}_{i_\alpha}\Bigl(w_h^{(k)}(\x^{(k,\ell,\alpha)})-w_h^{(\ell)}(\x^{(k,\ell,\alpha)})\Bigr).
    \end{align*}
    By assumption, $w_h$ satisfies the primal constraints and therefore we have $w_h^{(k)}(\x)=w_h^{(\ell)}(\x)$ for all vertices $\x \in \VV$. This concludes the proof.
\end{Proof}

With Lemma \ref{lem:w}, we can directly prove the conditions of Theorem \ref{thm:abstractcondition} in the next Lemma.

\begin{Lemma}\label{lem:winW}
    Let $D$ be the scaling matrix defined by selection scaling, $w_h = (w_h^{(1)\top},\ldots,w_h^{(K)\top})^\top \in \widetilde W_h$ and $v_h = (v_h^{(1)\top},\ldots,v_h^{(K)\top})^\top \in \widetilde W_h$ such that
    $
        \uv_h = DB_\Gamma^\top B_\Gamma \uw_h.
    $
    Then $v_h \in \widetilde W_{h,\Delta}$. If additionally $w_h \in \widetilde W_{h,\Delta}$, we further obtain that $B_\Gamma \uw_h = B_\Gamma \uv_h$.
\end{Lemma}
\begin{Proof}
    By Lemma \ref{lem:w}, we know that $v_h$ is the jump of $w_h$ across interfaces. Since $w_h \in \widetilde{W}_h$ is continuous at vertices, we know that $v_h$ must vanish at those vertices, i.e., $v_h \in \widetilde W_{h,\Delta}$, which proves the first part.
    
For the second part, let $w_{h} \in \widetilde W_{h,\Delta}$ and $v_h \in \widetilde W_{h,\Delta}$ such that $\uv_{h} = D B_\Gamma^\top B_\Gamma \uw_{h}$. Analogously to before, we have
    \begin{align*}
        v_i^{(k)} = \begin{cases}
            \lambda_{\Zc((k,i))} \quad &\mbox{ if } (k,i) \in \Ic_\Zc,\\
            0 \quad &\mbox{ else.}
        \end{cases}
    \end{align*}
    Suppose $r \notin \text{Im}(\Zc)$. Then $r \in \Rc_\Xc$ since $\{1,\ldots,n_\lambda\} = \Rc_\Xc \cup \text{Im}(\Zc)$. Note that if $r \in \Rc_\Xc$, we obtain $(B_\Gamma \uw_{h})_r = 0$ since $w_{h} \in \widetilde W_{h,\Delta}$. The first part proves that also $v_{h} \in W_{h,\Delta}$ and therefore we also obtain $(B_\Gamma \uv_{h})_r = 0$, which shows the desired result for $r \notin \text{Im}(\Zc)$. Consider now $r \in \text{Im}(\Zc)$. Applying $B_\Gamma$ to $\uv_{h}$ yields
    \begin{align*}
        (B_\Gamma \uv_{h,\Delta})_r = \sum_{k=1}^K \sum_{j \in \Ic_\Gamma^{(k)}} b_{r,j}^{(k)} v_j^{(k)}.
    \end{align*} 
    Summands only vanish if either $b_{r,j}^{(k)} = 0$ or $v_j^{(k)} = 0$. $v_j^{(k)}$ does not vanish if $j \in \Ic_\Zc$, that is, $\Rc_j^{(k)} = \{\Zc((k,j))\}$ and $b_{\Zc((k,j)),j}^{(k)} = 1$. However, this implies that $b_{r,j}^{(k)} = 0$ for $r \neq \Zc((k,j))$. Both conditions are true if $r \in \text{Im}(\Zc)$, in which case both summands do not vanish for the degree of freedom $j_0$ in patch $\Omega_{k_0}$ with $r = \Zc((k_0,j_0))$ and $b_{r,j_0}^{(k_0)} = 1$. Lemma \ref{lem:bijection} guarantees that $j_0$ is unique and, therefore, we obtain
    \begin{align*}
        (B_\Gamma \uv_{h,\Delta})_r = \underbrace{b_{r,j_0}^{(k_0)}}_{=1} v_{j_0}^{(k_0)} \lambda_r = (B_\Gamma \uw_{h,\Delta})_r
    \end{align*}
    for $r \in \text{Im}(\Zc)$, which concludes the proof.
\end{Proof}

Since both conditions of Theorem \ref{thm:abstractcondition} are satisfied, our next goal is to bound the supremum on the right-side of~\eqref{eq:omega}. In order to prove such a bound, we rely on results from \cite{Schneckenleitner2020}. The proof is divided into three parts: bounds for the discrete harmonic extension $\Hc_h^{(k)}$, an embedding result in the space $H^1(\Omega_k)$ and a tearing lemma in the space $H^{1/2}(\partial\Omega_k)$. We now state the main tools in order to prove the condition number estimate.
\begin{Lemma}\label{lem:boundsharmonic}
    Let $w_h^{(k)} \in W_h^{(k)}$ and $\Hc_h^{(k)} : W_h^{(k)} \to V_h^{(k)}$ be the discrete harmonic extension. Then
    \[
        |w_h^{(k)}|^2_{H^{1/2}(\partial \Omega_k)} \lesssim |\Hc_h^{(k)} w_h^{(k)}|^2_{H^1(\Omega_k)}\lesssim p\, |w_h^{(k)}|^2_{H^{1/2}(\partial \Omega_k)}.
    \]
\end{Lemma}
\begin{Proof}
    See \cite[Theorem 4.2]{Schneckenleitner2020}.
\end{Proof}
\begin{Lemma}[Embedding lemma]\label{lem:embedding}
    For $u \in V^{(k)}$ it holds that
    \[
        \inf_{c \in \RR} \sup_{x \in \cl\Omega_k} |u(x)-c|^2 \lesssim \left(1+\log p +\log \tfrac{H_k}{h_k}\right) |u|^2_{H^1(\Omega_k)}.
    \]
\end{Lemma}
\begin{Proof}
    The proof is done by using \cite[Lemma 4.14]{Schneckenleitner2020} and estimating the $L^2$ term by a Poincaré-type inequality.
\end{Proof}
The next lemma is an adapted version of the tearing lemma in \cite{Schneckenleitner2020}, which handles the conforming case where we estimate the $H^{1/2}$-norm on the whole boundary of the patch by tearing it apart into its $4$ edges. In the case of non-matching patches and the occurrence of T-junctions, one must take extra care. However, we only need the tearing lemma for functions that are non-zero on the refined sides of interfaces, which are again full edges of a patch and can be handled by the already established results. We formulate our tearing lemma in the following lemma:
\begin{Lemma}[Tearing lemma]\label{lem:tearing}
 Let $v_h \in V_h^{(k)}$ such that $\ell \nprec k  \implies v_h|_{\Gamma^{(k,\ell)}} = 0  $. Then the following estimates hold:
\begin{align*}
    &|v_h|_{H^{1/2}(\partial \Omega_k)}^2 \lesssim \sum_{\ell \in N_\Gamma(k)} \left(|v_h|_{H^{1/2}(\Gamma^{(k,\ell)})}^2 + \Lambda^{(k)}| v_h|_{L_\infty^0(\Gamma^{(k,\ell)})}^2\right),
\end{align*}
where $\Lambda^{(k)} \coloneqq 1+\log p +\log \tfrac{H_k}{h_k}$ and $|v_h|_{L^0_\infty(T)} \coloneqq \inf_{c \in \RR} \sup_{x \in T}|v_h(x)-c|$. Moreover, for all $v_h \in V_h^{(k)}$, the estimate
\[
    \sum_{\ell \in N_\Gamma(k)} |v_h|_{H^{1/2}(\Gamma^{(k,\ell)})}^2 
    \le
    |v_h|_{H^{1/2}(\partial \Omega_k)}^2
\]
holds.
\end{Lemma}
\begin{Proof}
    See \cite[Lemma 4.15]{Schneckenleitner2020} since $\Gamma^{(k,\ell)}$ is a full edge of $\Omega_k$ (i.e., the image of one of the four sides of the unit square $\widehat\Omega = (0,1)^2$ under $G_k$) if $\ell \prec k$.
\end{Proof}

Now, we state the main result of this paper.
\begin{Theorem}\label{thrm:final}
    Let $F$ and $M_{sD}$ be given as in \eqref{eq:F} and \eqref{eq:MsD}, where $D$ is given by selection scaling~\eqref{eq:D def}. Then, there is a generic constant $C$ that 
    only depends
    on the constant $C_G$ from Assumption~\ref{Ass:GeoEquiv}
    and
    the constant $C_Q$ from Assumption~\ref{Ass:QuasiUniform}
    such that the condition number of the preconditioned system satisfies 
    \begin{equation}\label{eq:coefficientdependent}
        \kappa(M_{sD}F) \le C \, p \max\left\{\max_{\ell \prec k} \tfrac{\nu^{(k)}}{\nu^{(\ell)}},1\right\} \left(1+\log p + \max_{k=1,\ldots,K} \log \frac{H_k}{h_k}\right)^2.
    \end{equation}
\end{Theorem}
\begin{Proof}
Let $w_h=(w_h^{(1)},\ldots,w_h^{(K)}) \in \widetilde W_h$ with coefficient vector $\uw_h $ and $v_h = (v_h^{(1)},\ldots,v_h^{(K)})$ with coefficient vector $\uv_h$ such that $\uv_h = D B_\Gamma^\top B_\Gamma \uw_h$. By Theorem \ref{lem:boundsharmonic} we immediately obtain
\begin{align*}
    \| D B_\Gamma^\top B \uw_h \|_S ^2
        & = \|\uv_h\|_S ^2 = \sum_{k=1}^K \| \uv_h^{(k)}\|_{S^{(k)}} ^2 = \sum_{k=1}^K \nu^{(k)}|\Hc_h^{(k)} v_h^{(k)}|_{H^1(\Omega_k)}^2 \lesssim p\sum_{k=1}^K \nu^{(k)} |v_h^{(k)}|_{H^{1/2}(\partial \Omega_k)}^2,
\end{align*}
where $\Hc_h^{(k)}$ is the discrete harmonic extension into $H^1(\Omega_k)$. We define $\Lambda \coloneqq \max_{1,\ldots,K}\Lambda^{(k)}$ and by using Lemma \ref{lem:tearing} and Lemma \ref{lem:w}, we further get
\begin{align*}
       \| D B_\Gamma^\top B \uw_h \|_S ^2 &\lesssim p \Lambda \sum_{k=1}^K \nu^{(k)} \sum_{\ell \prec k} \left(|v_h^{(k)}|_{H^{1/2}(\Gamma^{(k,\ell)})}^2 + |v_h^{(k)}|_{L_\infty^0(\Gamma^{(k,\ell)})}^2\right) \\
        &\lesssim p \Lambda \sum_{k=1}^K \nu^{(k)} \sum_{\ell \prec k} \left(|w_h^{(k)}-w_h^{(\ell)}|_{H^{1/2}(\Gamma^{(k,\ell)})}^2 + |w_h^{(k)}-w_h^{(\ell)}|_{L_\infty^0(\Gamma^{(k,\ell)})}^2\right) \\      
        &\lesssim  p C_\nu \Lambda \sum_{k=1}^K \nu^{(k)}\left(|w_h^{(k)}|_{H^{1/2}(\partial \Omega_k)}^2 + |w_h^{(k)}|_{L_\infty^0(\partial \Omega_k)}^2\right),
\end{align*}
where $C_\nu \coloneqq \max\left\{\max_{\ell \prec k} \tfrac{\nu^{(k)}}{\nu^{(\ell)}},1\right\}$. Applying again Lemma \ref{lem:boundsharmonic} and using the definition of the seminorm $|\cdot|_{L_\infty^0(\partial \Omega_k)}$, we obtain
\begin{align*}
        \| D B_\Gamma^\top B \uw_h \|_S ^2 &\lesssim p C_\nu \Lambda \sum_{k=1}^K \nu^{(k)}\left(|\Hc_h^{(k)}w_h^{(k)}|_{H^{1}(\Omega_k)}^2 + \inf_{c \in \RR} \sup_{x \in \partial\Omega_k}|w_h^{(k)}(x)-c|^2\right).
\end{align*}
Using \eqref{eq:harmonicextension} and since $\partial\Omega_k \subset \cl\Omega_k$ it also holds that
\[
    \sup_{x \in \partial\Omega_k}|w_h^{(k)}(x)| \leq \sup_{x \in \cl\Omega_k}|\Hc_h^{(k)}w_h^{(k)}(x)| \quad \mbox{for all } w_h^{(k)} \in W_h^{(k)}.
\]
Therefore, we further use Lemma \ref{lem:embedding}, the linearity of $\Hc_h^{(k)}$, $\Hc_h^{(k)}c = c$ and the embedding Lemma \ref{lem:embedding} and obtain
\begin{align*}
\| D B_\Gamma^\top B \uw_h \|_S ^2 & \lesssim p C_\nu \Lambda \sum_{k=1}^K \nu^{(k)}\left(|\Hc_h^{(k)}w_h^{(k)}|_{H^{1}(\Omega_k)}^2 + \inf_{c \in \RR} \sup_{x \in \cl\Omega_k} |(\Hc_h^{(k)}w_h^{(k)})(x)-c|^2\right)\\
         & \lesssim p C_\nu \Lambda \sum_{k=1}^K \nu^{(k)}\left((1+2\Lambda)|\Hc_h^{(k)}w_h^{(k)}|_{H^{1}(\Omega_k)}^2\right)\\
         & \lesssim p C_\nu \Lambda^2 \sum_{k=1}^K \nu^{(k)}|\Hc_h^{(k)}w_h^{(k)}|_{H^{1}(\Omega_k)}^2 = p C_\nu \Lambda^2 \|\uw\|_S ^2.
\end{align*}
Using Theorem \ref{thm:abstractcondition}, we arrive at the desired result.
\end{Proof}

\begin{Remark}
    We that the condition number bound is not robust with respect to the diffusion coefficient $\nu$, specifically, the ratio of $\nu$ across different patches. The numerical tests confirm this. The problem is that the scaling matrix $D$ only considers a specific side of an interface $\Gamma^{(k,\ell)}$, i.e $F^{-1} \approx \sum_{\ell \prec k} B_{\Gamma^{(k,\ell)}}^{(k)} S_{\Gamma^{(k,\ell)}}^{(k)} B_{\Gamma^{(k,\ell)}}^{(k)\top}$, where $\square_{\Gamma^{(k,\ell)}}^{(k)}$ means, we only take degrees of freedom in $\Omega_k$ on the interface $\Gamma^{(k,\ell)}$. In the case of matching interfaces, one may circumvent this problem by employing \textit{coefficient scaling} as a suitable scaling technique.
    Alternatively, one can define the relation $\prec$ with the coefficient in mind. However, these techniques cannot be applied in a straightforward way to non-matching interfaces.
    A remedy to this problem in the case of non-matching interfaces would be a deluxe type preconditoner that would somehow involve both sides of an interface, i.e., 
    \begin{equation}\label{eq:edgeprec} F^{-1}\approx \sum_{\ell \prec k} B_{\Gamma^{(k,\ell)}}^{(k)}\left(S_{\Gamma^{(k,\ell)}}^{(k)-1} + P_{\ell}^{k} S_{\Gamma^{(k,\ell)}}^{(\ell)-1} P_{k}^{\ell}\right)^{-1} B_{\Gamma^{(k,\ell)}}^{(k)\top},
    \end{equation} where $P_{\ell}^{k}$ is the prolongation matrix for the canonical embedding and $P_{k}^{\ell}:=(P^{k}_{\ell})^\top$.
\end{Remark}

\section{Numerical results}\label{sec:num}

In this section, we present numerical tests illustrating the theory developed in the previous section.

\subsection{Checkerboard pattern}
\begin{figure}[ht]\label{fig:checker_uni_good}
    \begin{subfigure}[t]{0.3\textwidth}
        \centering
        \includegraphics[width=0.8\textwidth]{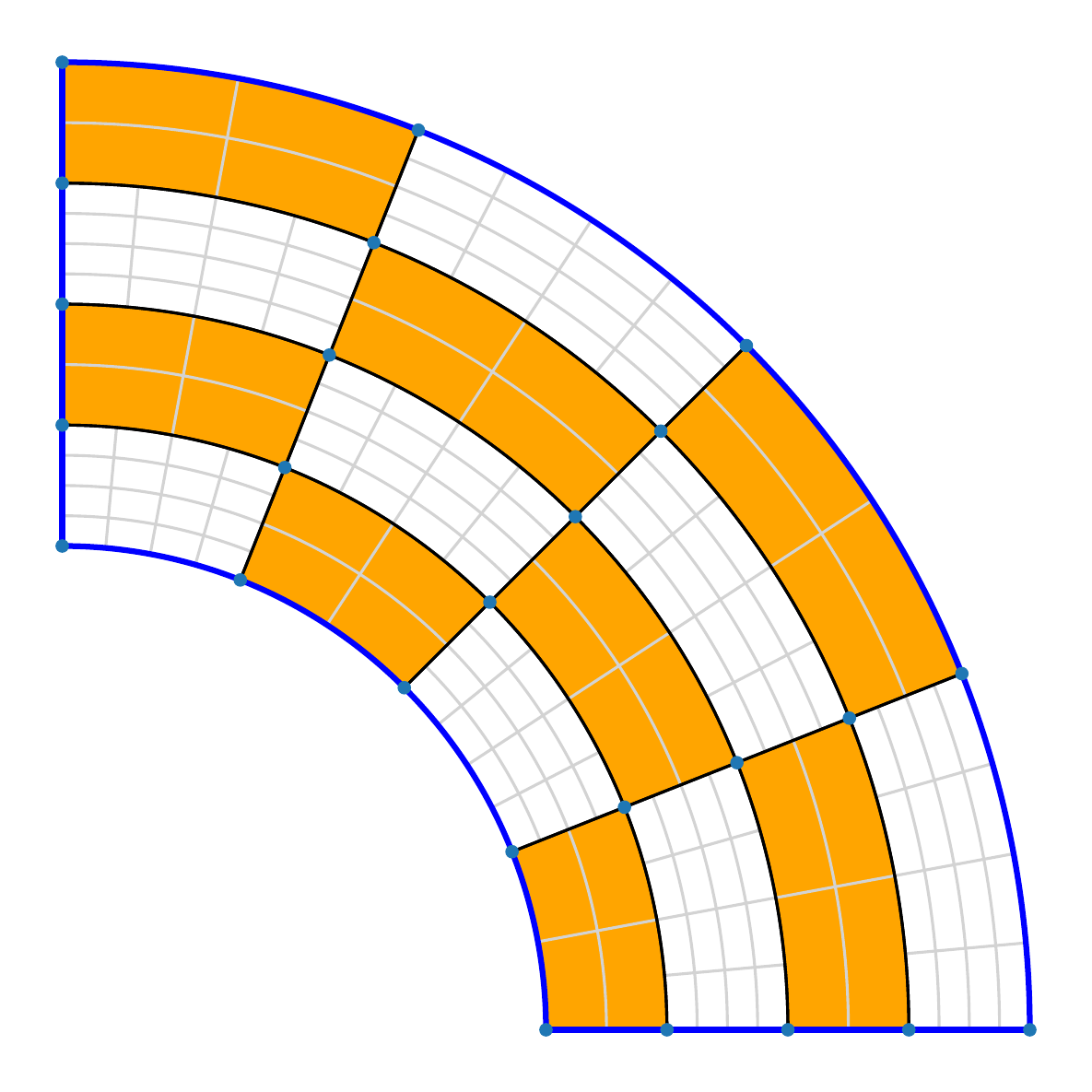}
        \caption{$r=0$}
    \end{subfigure}
    \hfill
    \begin{subfigure}[t]{0.3\textwidth}
        \centering
        \includegraphics[width=0.8\linewidth]{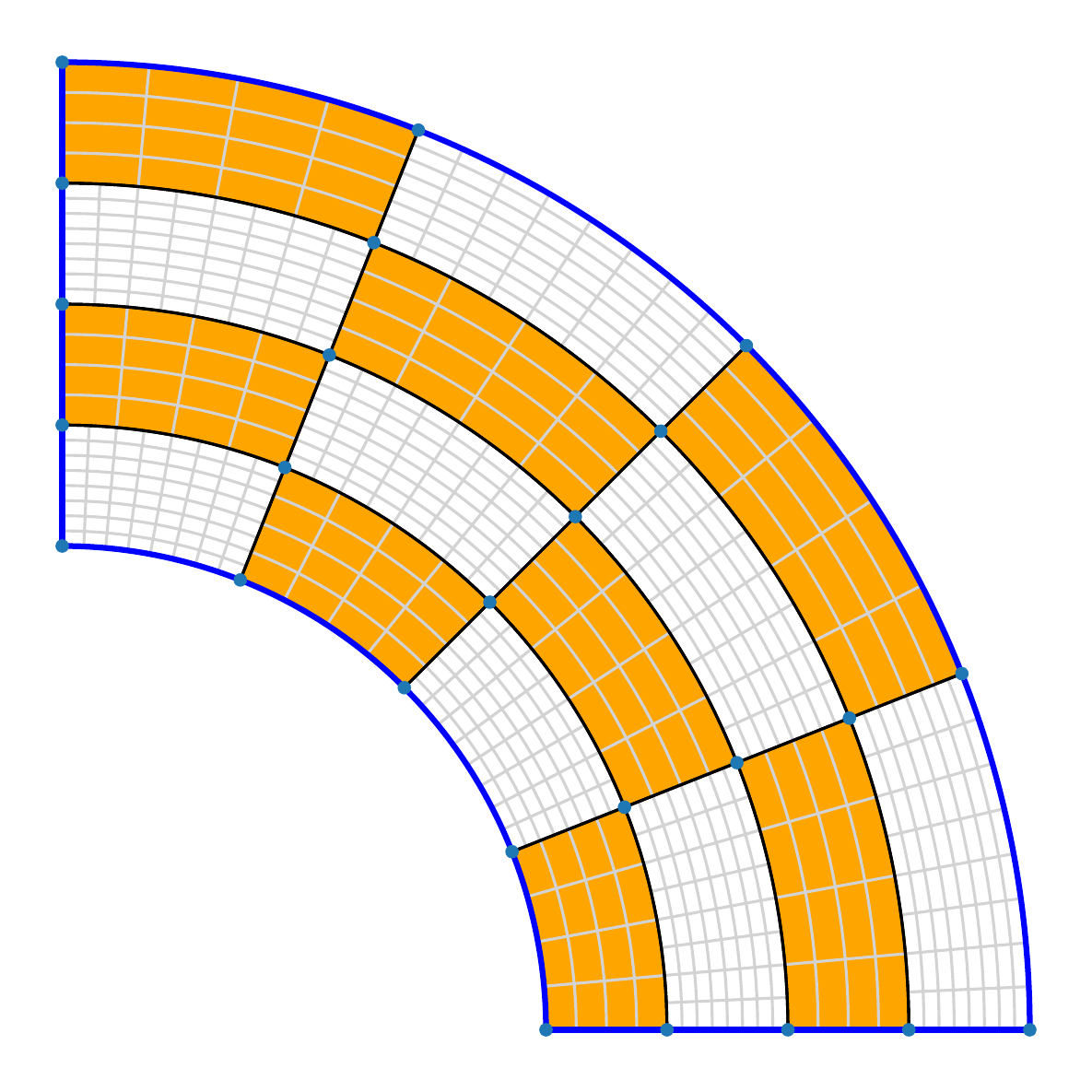}
        \caption{$r=1$}
    \end{subfigure}
    \hfill
    \begin{subfigure}[t]{0.3\textwidth}
        \centering
        \includegraphics[width=0.8\linewidth]{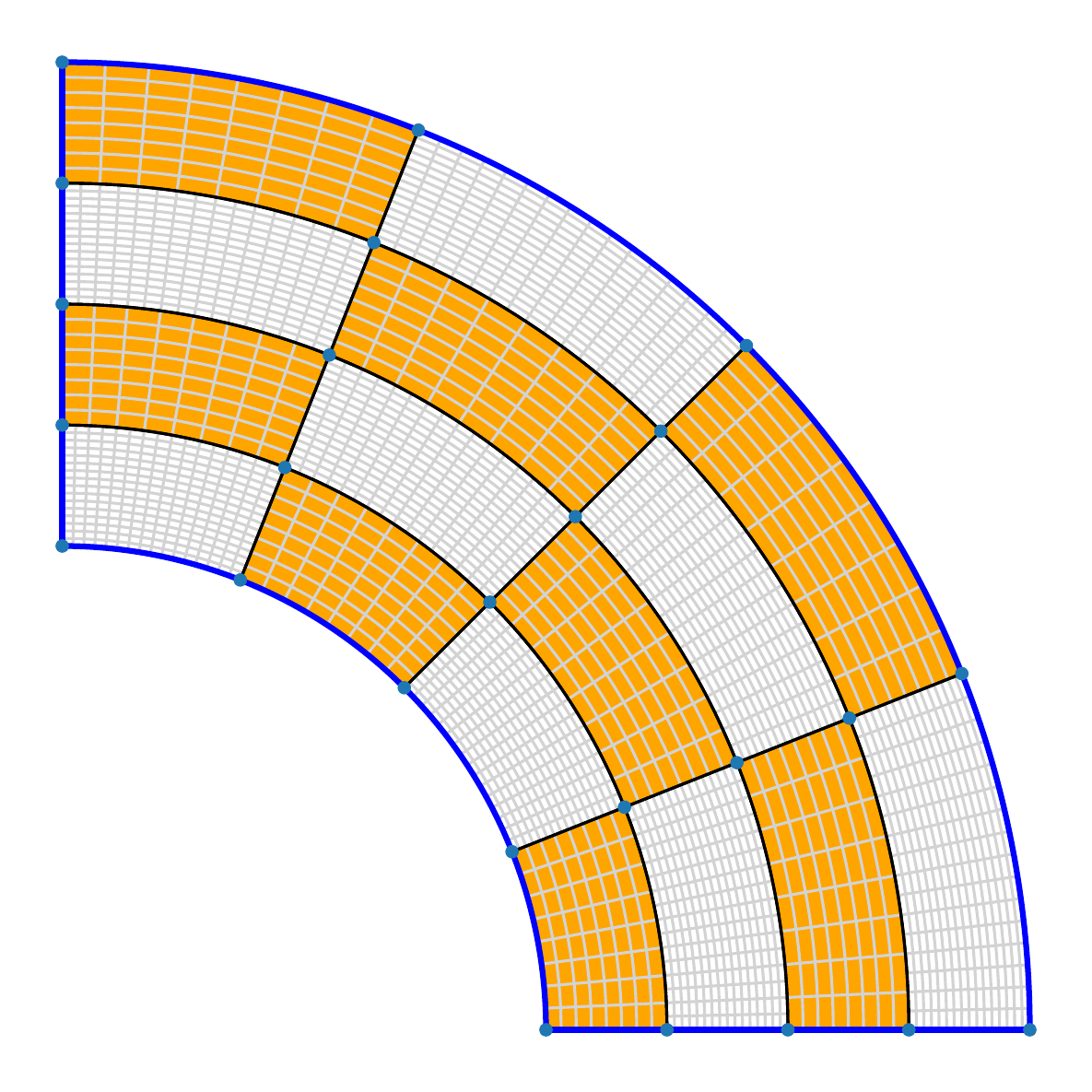}
        \caption{$r = 2$}
    \end{subfigure}
    \caption{Checkerboard multi-patch domains with mesh disparity level $d = 1$.}
    \label{fig:checkerboard}
\end{figure}

As a first test case, we set up a simple multi-patch domain without any T-junctions. Specifically, the computational domain $\Omega$ is a quarter annulus, composed of $16$ patches in a checkerboard style, as shown in Figure~\ref{fig:checkerboard}. 
We consider the following model problem: Find $u \in H_0^{1}(\Omega)$ such that
\begin{equation}\label{eq:nu poisson}
\int_\Omega \nu \nabla u \cdot \nabla v \diff x = \int_\Omega f v \diff x \quad \text{for all } v \in H_0^1(\Omega)
\quad \mbox{with}\quad f=1.
\end{equation}
We choose $\nu$ to be patchwise constant. Specifically, we choose $\nu_{\text{orange}} = 1000$ for the patches depicted in orange color and $\nu_{\text{white}}=1$ for the patches depicted in white color. 
For the discretization, we choose the initial grid sizes $\wh h_0 = (p+1)^{-1}$ and apply $r$ levels of uniform refinement. Some of the patches are refined more often (as will be discussed below). The maximum difference of the refinement levels between neighboring patches is called the mesh level disparity $d$. 
For all experiments, we present the number of iterations ("it") that the conjugate gradient solver required to reduce the residual by a factor of $10^{-6}$. Moreover, the condition number of the preconditioned system ("$\kappa$"), as estimates by the conjugate gradient solver, are presented. In all cases, we apply the conjugate gradient solver to the Schur complement formulation~\eqref{eq:schur}.

For the first experiments, we refine the patches depicted in white color once more ($d=1$), see Figure~\ref{fig:checkerboard}.
In Table~\ref{tab:checkerboard:selection}, we provide the results for the scaled Dirichlet preconditioner with selection scaling. Table~\ref{tab:checkerboard:edge} shows the analogous results if the deluxe type preconditioner~\eqref{eq:edgeprec} is applied. Note that the diffusion parameters are chosen in a "good" way, since in this case we have 
$$\max\left\{\max_{\ell \prec k} \tfrac{\nu^{(k)}}{\nu^{(\ell)}},1\right\} = 1$$
by construction, which is the coefficient dependent term in \eqref{eq:coefficientdependent}. In Table \ref{tab:checkerboard:selection}, we can actually see that selection scaling is near-optimal in this case and the results are surprisingly even better than what \eqref{eq:coefficientdependent} suggests.

\begin{table}[ht]
    \centering
    \begin{tabular}{c|lr|lr|lr|lr|lr|lr|lr}
	\toprule
	& \multicolumn{2}{c|}{$p=2$}
	& \multicolumn{2}{c|}{$p=3$}
	& \multicolumn{2}{c|}{$p=4$}
        & \multicolumn{2}{c|}{$p=5$}
        & \multicolumn{2}{c|}{$p=6$}
        & \multicolumn{2}{c|}{$p=7$}
	& \multicolumn{2}{c }{$p=8$}\\
	$\wh h$
	& it & $\kappa$
	& it & $\kappa$
	& it & $\kappa$
        & it & $\kappa$
        & it & $\kappa$
        & it & $\kappa$
	& it & $\kappa$\\
	\midrule
        $2^{-1}$ & 3  & 1.009 & 3  & 1.013 & 3 & 1.017 & 3 & 1.020 & 3 & 1.023 & 3 & 1.026 & 3 & 1.028\\
        $2^{-2}$ & 3 & 1.013 & 3 & 1.018 & 3 & 1.022 & 3 & 1.025 & 3 & 1.029 & 3 & 1.032 & 3 & 1.034\\
        $2^{-3}$ & 3 & 1.017 & 3 & 1.023 & 3 & 1.027 & 3 & 1.031 & 3 & 1.035 & 3 & 1.038 & 3 & 1.041\\
        $2^{-4}$ & 3 & 1.022 & 3 & 1.028 & 3 & 1.034 & 3 & 1.038 & 3 & 1.042 & 3 & 1.046 & 3 & 1.049\\
        $2^{-5}$ & 3 & 1.028 & 3 & 1.035 & 3 & 1.041 & 3 & 1.045 & 3 & 1.050 & 3 & 1.054 & 3 & 1.057\\
        $2^{-6}$ & 3 & 1.034 & 3 & 1.042 & 3 & 1.048 & 3 & 1.053 & 3 & 1.058 & \multicolumn{2}{c|}{OoM} & \multicolumn{2}{c}{OoM}\\
	\bottomrule
    \end{tabular}
    \caption{\label{tab:checkerboard:selection}Iteration counts (it) and condition numbers $\kappa$ for the checkerboard mesh with scaled Dirichlet preconditioner with selection scaling and mesh level disparity $d = 1$.}
\end{table}

\begin{table}[ht]
    \centering
    \begin{tabular}{c|lr|lr|lr|lr|lr|lr|lr}
	\toprule
	& \multicolumn{2}{c|}{$p=2$}
	& \multicolumn{2}{c|}{$p=3$}
	& \multicolumn{2}{c|}{$p=4$}
        & \multicolumn{2}{c|}{$p=5$}
        & \multicolumn{2}{c|}{$p=6$}
        & \multicolumn{2}{c|}{$p=7$}
	& \multicolumn{2}{c }{$p=8$}\\
	$\wh h$
	& it & $\kappa$
	& it & $\kappa$
	& it & $\kappa$
        & it & $\kappa$
        & it & $\kappa$
        & it & $\kappa$
	& it & $\kappa$\\
	\midrule
        $2^{-1}$ & 15 & 10.77 & 17 & 15.33 & 18 & 19.23 & 19 & 22.72 & 19 & 25.78 & 20 & 28.62 & 21 & 31.17\\
        $2^{-2}$ & 17 & 14.93 & 18 & 20.16 & 19 & 24.61 & 20 & 28.48 & 21 & 31.90 & 23 & 35.01 & 23 & 37.73\\
        $2^{-3}$ & 18 & 19.71 & 19 & 25.64 & 21 & 30.60 & 23 & 34.86 & 23 & 38.49 & 24 & 41.88 & 24 & 44.77\\
        $2^{-4}$ & 19 & 25.12 & 22 & 31.74 & 23 & 37.11 & 24 & 41.73 & 25 & 45.69 & 27 & 49.40 & 27 & 52.47\\
        $2^{-5}$ & 22 & 31.17 & 23 & 38.31 & 24 & 44.11 & 27 & 49.23 & 27 & 53.33 & 27 & 56.90 & 28 & 60.46\\
        $2^{-6}$ & 23 & 37.70 & 26 & 45.51 & 27 & 51.77 & 27 & 56.73 & \multicolumn{2}{c|}{OoM} & \multicolumn{2}{c|}{OoM} & \multicolumn{2}{c}{OoM}\\
    \bottomrule
    \end{tabular}
    \caption{\label{tab:checkerboard:edge}Iteration counts (it) and condition numbers $\kappa$ for the checkerboard mesh with deluxe type preconditioner and mesh level disparity $d = 1$.}
\end{table}

\begin{figure}[ht]
    \begin{subfigure}[t]{0.3\textwidth}
        \centering
        \includegraphics[width=0.8\textwidth]{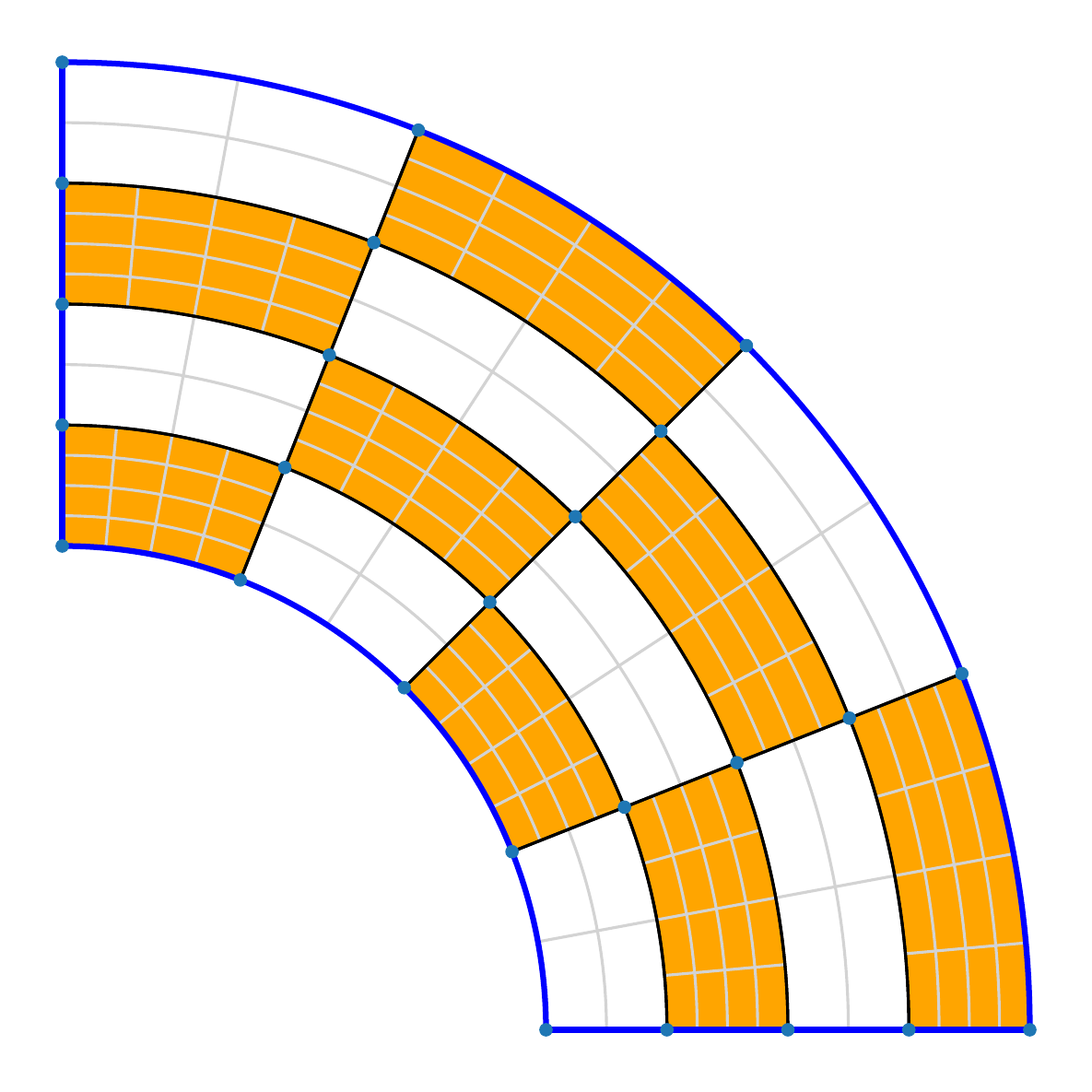}
        \caption{$r=0$}
    \end{subfigure}
    \hfill
    \begin{subfigure}[t]{0.3\textwidth}
        \centering
        \includegraphics[width=0.8\linewidth]{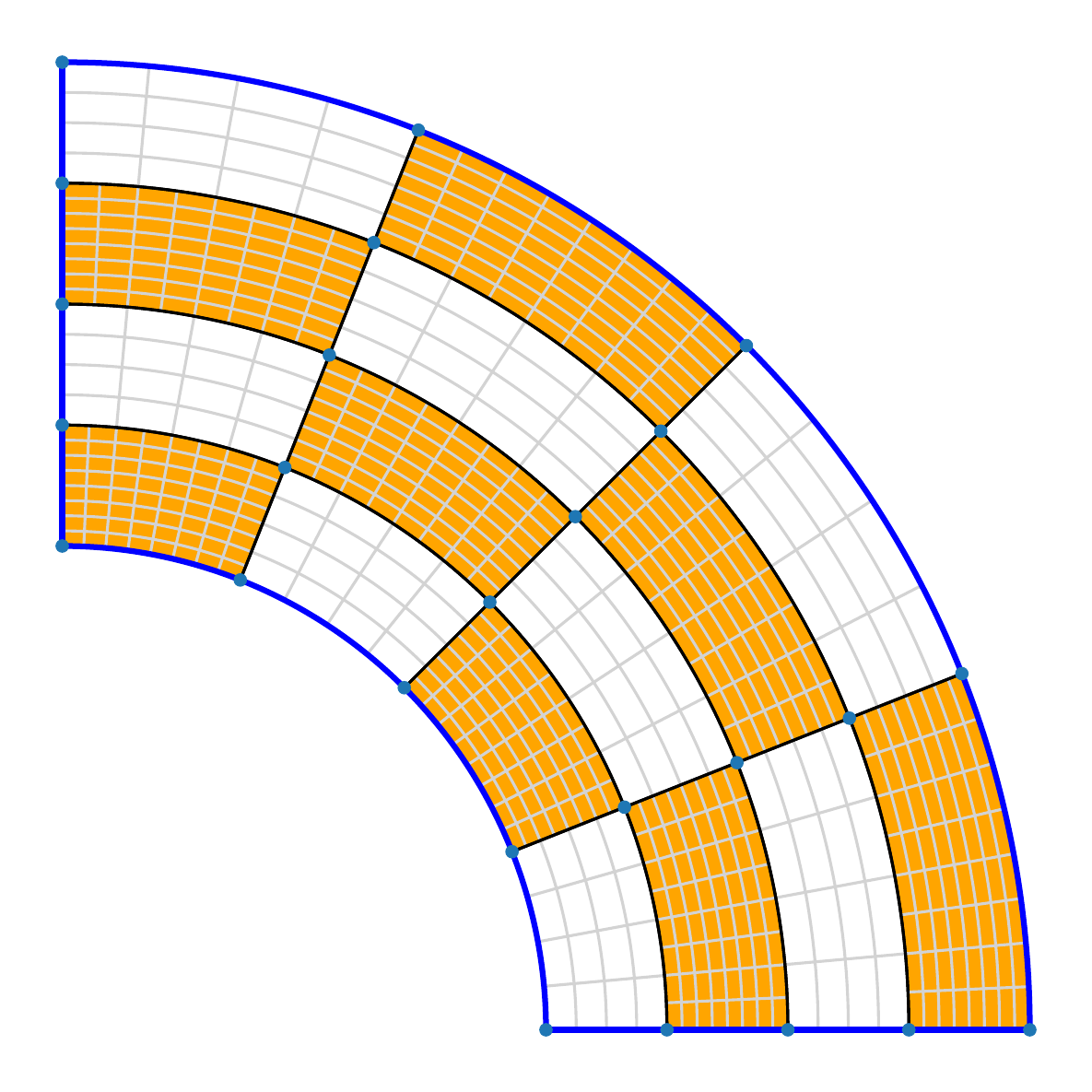}
        \caption{$r=1$}
    \end{subfigure}
    \hfill
    \begin{subfigure}[t]{0.3\textwidth}
        \centering
        \includegraphics[width=0.8\linewidth]{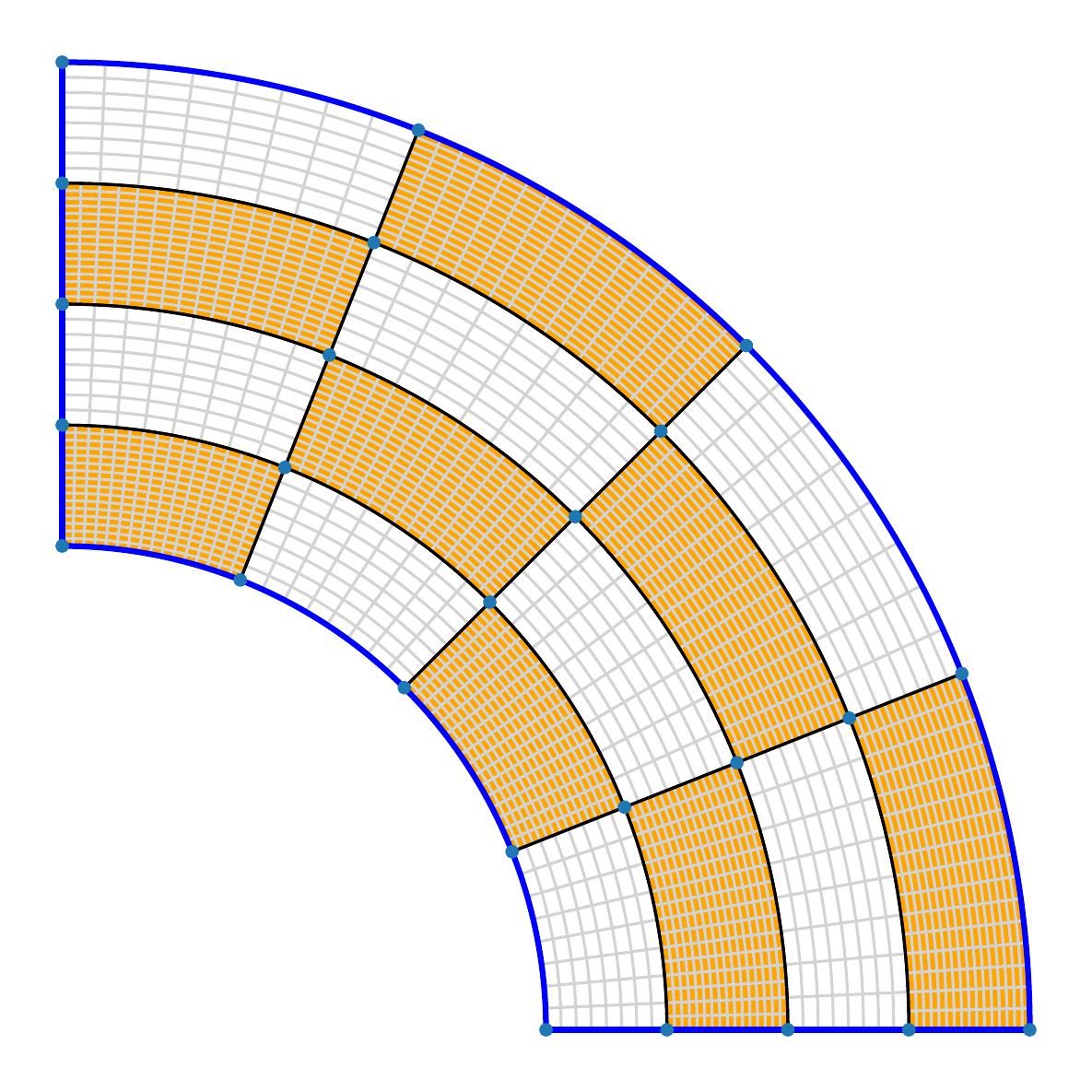}
        \caption{$r=2$}
    \end{subfigure}
    \caption{Checkerboard multi-patch domains with mesh disparity level $d = 1$.}
    \label{fig:checkerboard2}
\end{figure}

For the next experiments, we still have $\nu_{\text{orange}} = 1000$, $\nu_{\text{white}} = 1$. As opposed to the first experiments, the grid used on the orange-colored patches is refined once more than the grid on the white-colored patches, see Figure~\ref{fig:checkerboard2}. This is the "bad" case, where
$$\max\left\{\max_{\ell \prec k} \tfrac{\nu^{(k)}}{\nu^{(\ell)}},1\right\} = \max_{\ell \prec k} \tfrac{\nu^{(k)}}{\nu^{(\ell)}} = 1000.$$
The iteration numbers and condition numbers are displayed in Table~\ref{tab:checkerboard2:selection}, from which we observe that the condition numbers are increased by a factor of $1000$. Furthermore, we note that the deluxe type preconditioner \eqref{eq:edgeprec} is consistent in both examples (see Tables \ref{tab:checkerboard:edge} and \ref{tab:checkerboard2:edge}), highlighting its robustness with respect to the diffusion parameter. 

\newlength{\digitwidth}
\settowidth{\digitwidth}{$0$}

\begin{table}[ht]
    \centering
    \setlength{\tabcolsep}{5pt}
    \begin{tabular}{c|lr|lr|lr|lr|lr|lr|lr}
	\toprule
	& \multicolumn{2}{c|}{$p=2$}
	& \multicolumn{2}{c|}{$p=3$}
	& \multicolumn{2}{c|}{$p=4$}
        & \multicolumn{2}{c|}{$p=5$}
        & \multicolumn{2}{c|}{$p=6$}
        & \multicolumn{2}{c|}{$p=7$}
	& \multicolumn{2}{c }{$p=8$}\\
	$\wh h$
	& it & $\kappa$
	& it & $\kappa$
	& it & $\kappa$
        & it & $\kappa$
        & it & $\kappa$
        & it & $\kappa$
	& it & $\kappa$\\
	\midrule
$2^{-1}$ & 49 & 6.0e3 & 53 & 8.7e3 & 63 & 1.1e4 & 70 & 1.3e4 & 77 & 1.5e4 & 83 & 1.7e4 & 90 & 1.9e4\\
$2^{-2}$ & 61 & 8.4e3 & 66 & 1.2e4 & 74 & 1.4e4 & 84 & 1.7e4 & 92 & 1.9e4 & 98 & 2.1e4 & 104 & 2.3e4\\
$2^{-3}$ & 64 & 1.1e4 & 77 & 1.5e4 & 88 & 1.8e4 & 97 & 2.1e4 & 104 & 2.4e4 & 108 & 2.6e4 & 114 & 2.8e4\\
$2^{-4}$ & 74 & 1.5e4 & 89 & 1.9e4 & 101 & 2.3e4 & 106 & 2.6e4 & 113 & 2.8e4 & 122 & 3.1e4 & 129 & 3.3e4\\
$2^{-5}$ & 88 & 1.9e4 & 104 & 2.3e4 & 112 & 2.7e4 & 121 & 3.1e4 & 130 & 3.4e4 & 135 & 3.6e4 & 144 & 3.9e4\\
$2^{-6}$ & 101 & 2.3e4 & 114 & 2.8e4 & 127 & 3.3e4 & 134 & 3.6e4 & \multicolumn{2}{c|}{OoM} & \multicolumn{2}{c|}{OoM} & \multicolumn{2}{c}{OoM}\\
	\bottomrule
    \end{tabular}
    \caption{\label{tab:checkerboard2:selection}Iteration counts (it) and condition numbers $\kappa$ for the checkerboard mesh with scaled Dirichlet preconditioner with selection scaling and mesh level disparity $d = 1$.}
\end{table}

\begin{table}[ht]
    \centering
    \begin{tabular}{c|lr|lr|lr|lr|lr|lr|lr}
	\toprule
	& \multicolumn{2}{c|}{$p=2$}
	& \multicolumn{2}{c|}{$p=3$}
	& \multicolumn{2}{c|}{$p=4$}
        & \multicolumn{2}{c|}{$p=5$}
        & \multicolumn{2}{c|}{$p=6$}
        & \multicolumn{2}{c|}{$p=7$}
	& \multicolumn{2}{c }{$p=8$}\\
	$\wh h$
	& it & $\kappa$
	& it & $\kappa$
	& it & $\kappa$
        & it & $\kappa$
        & it & $\kappa$
        & it & $\kappa$
	& it & $\kappa$\\
	\midrule
        $2^{-1}$ & 13 & 7.78 & 15 & 11.32 & 17 & 14.68 & 17 & 17.68 & 18 & 20.42 & 19 & 22.94 & 19 & 25.23\\
        $2^{-2}$ & 15 & 10.78 & 17 & 15.34 & 18 & 19.23 & 19 & 22.72 & 19 & 25.77 & 21 & 28.60 & 22 & 31.15\\
        $2^{-3}$ & 17 & 14.93 & 18 & 20.17 & 19 & 24.60 & 21 & 28.47 & 22 & 31.88 & 23 & 34.96 & 23 & 37.67\\
        $2^{-4}$ & 18 & 19.71 & 20 & 25.63 & 21 & 30.54 & 23 & 34.82 & 23 & 38.43 & 24 & 41.80 & 25 & 44.70\\
        $2^{-5}$ & 19 & 25.12 & 22 & 31.71 & 24 & 37.08 & 25 & 41.68 & 25 & 45.50 & 26 & 48.99 & 27 & 52.34\\
        $2^{-6}$ & 23 & 31.16 & 24 & 38.28 & 25 & 44.05 & 26 & 48.84 & \multicolumn{2}{c|}{OoM} & \multicolumn{2}{c|}{OoM} & \multicolumn{2}{c}{OoM}\\
	\bottomrule
    \end{tabular}
    \caption{Iteration counts (it) and condition numbers $\kappa$ for the checkerboard mesh with deluxe type preconditioning and mesh level disparity $d = 1$.}
    \label{tab:checkerboard2:edge}
\end{table}

\begin{figure}[ht]
    \begin{subfigure}[t]{0.3\textwidth}
        \centering
        \includegraphics[width=0.8\textwidth]{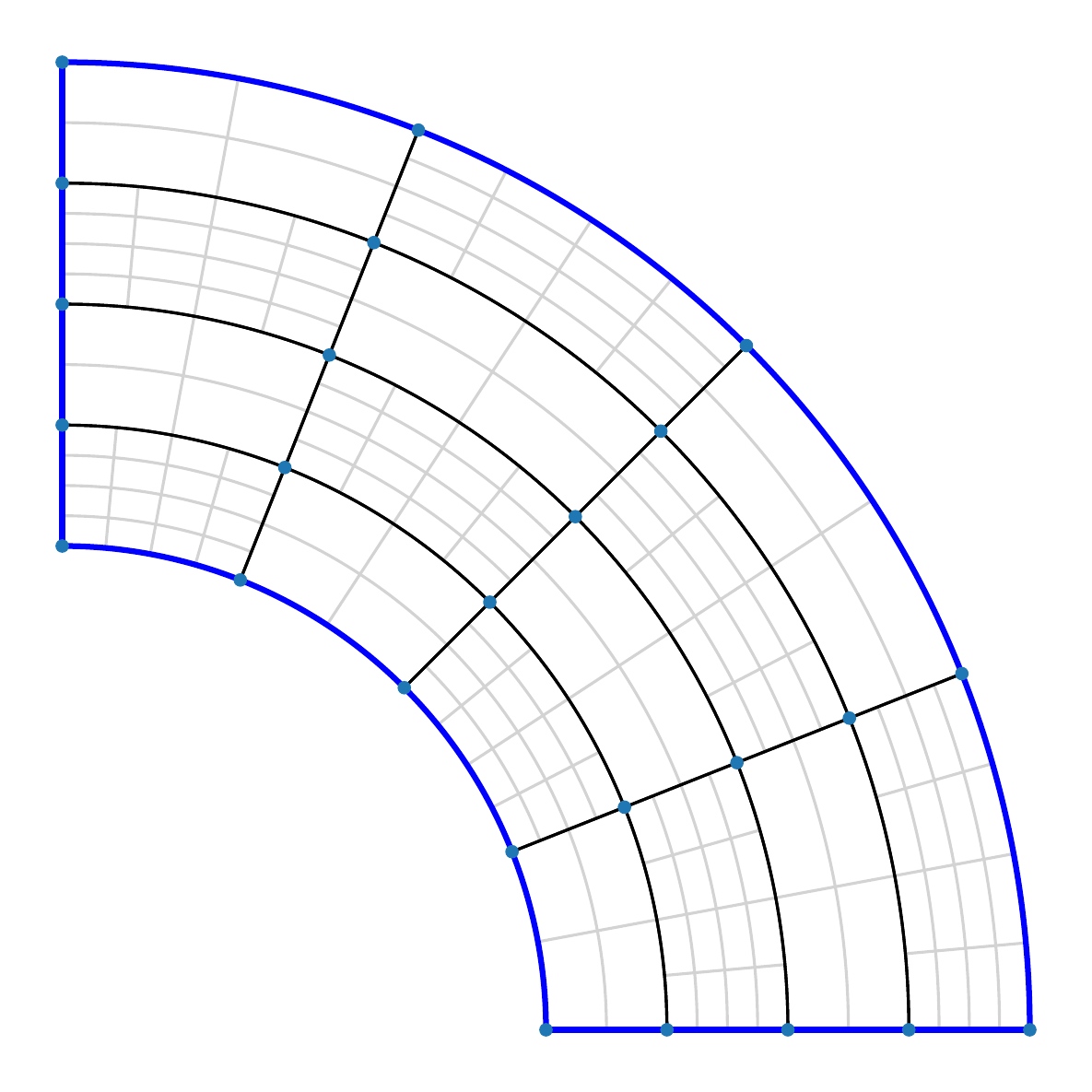}
        \caption{$d=1$}
    \end{subfigure}
    \hfill
    \begin{subfigure}[t]{0.3\textwidth}
        \centering
        \includegraphics[width=0.8\linewidth]{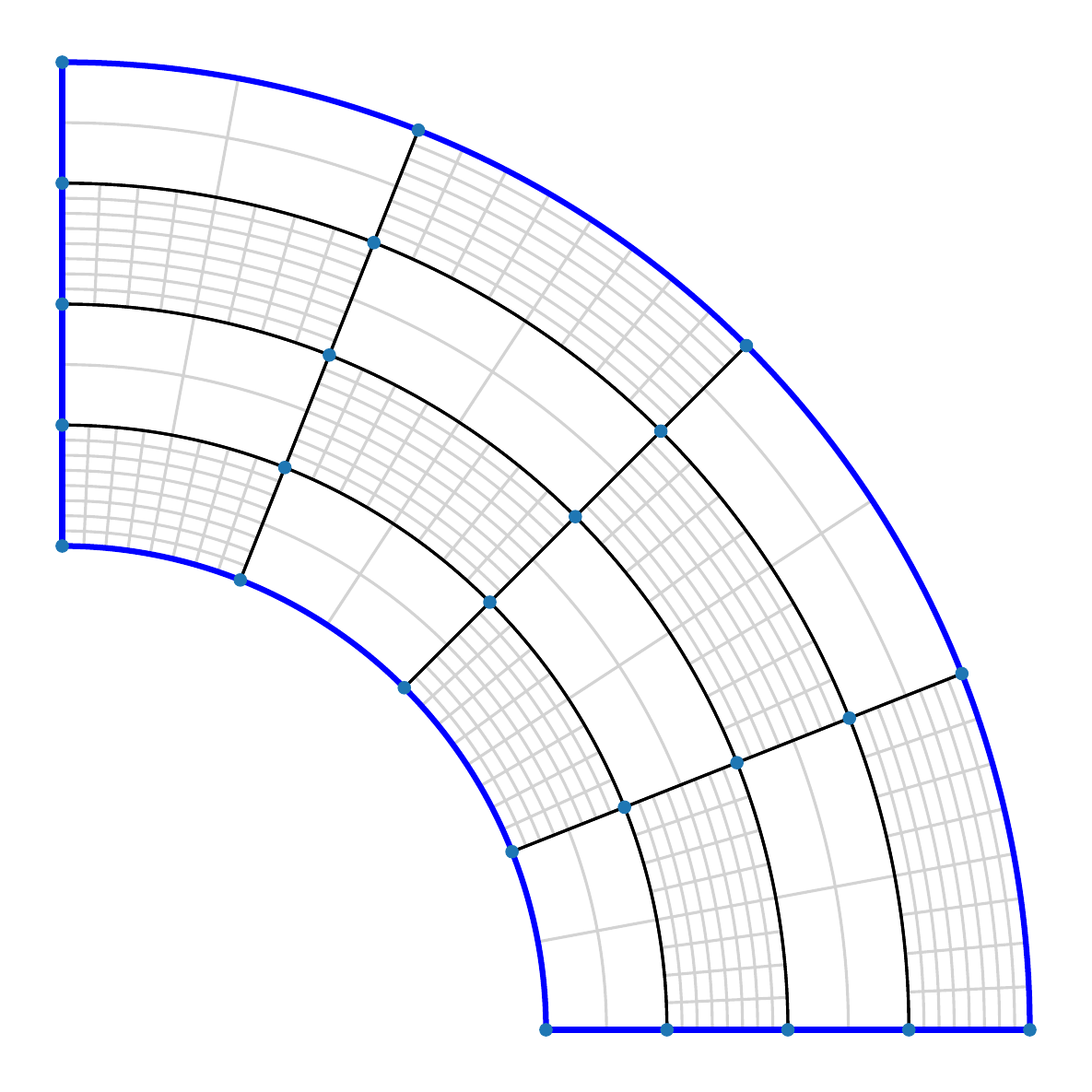}
        \caption{$d=2$}
    \end{subfigure}
    \hfill
    \begin{subfigure}[t]{0.3\textwidth}
        \centering
        \includegraphics[width=0.8\linewidth]{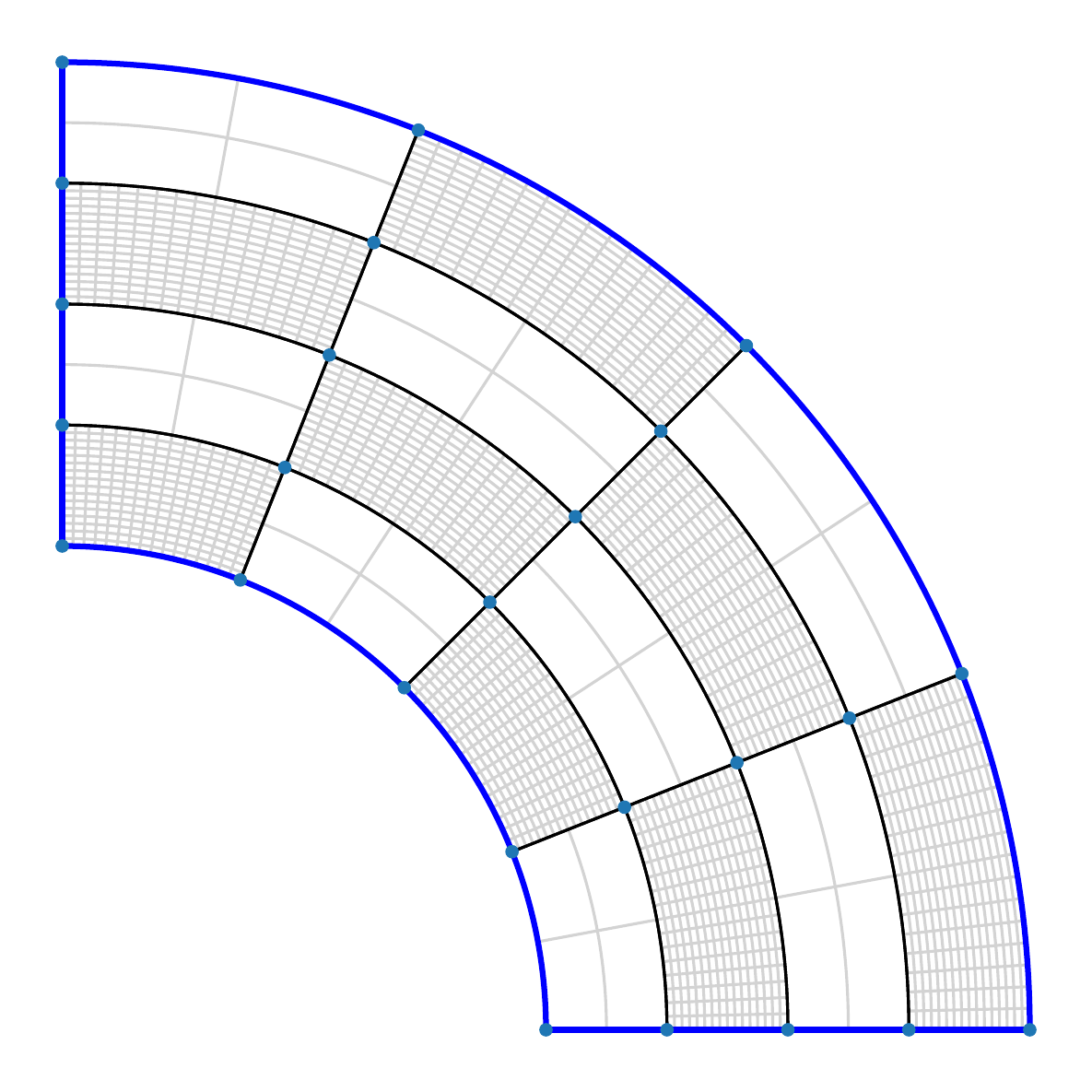}
        \caption{$d=3$}
    \end{subfigure}
    \caption{Checkerboard multi-patch domains with varying mesh disparity level $d$.}
    \label{fig:checker_dis}
\end{figure}

We also inspect the behavior of the solver if the mesh disparity level increases, i.e., we refine one set of patches in the checkerboard but leave the rest unchanged, as displayed in Figure \ref{fig:checker_dis}. Here, we set the diffusion parameter $\nu=1$ for the whole domain, since it influences the condition number only by a constant. In Table \ref{tab:checker_dis:selection}, we present the iteration number and condition number for the scaled Dirichlet preconditioner with selection scaling, and in Table \ref{tab:checker_dis:edge}, the same for the deluxe-type preconditioner.  We observe robustness in the case of the scaled Dirichlet preconditioner with selection scaling (as predicted by the theory), the deluxe type preconditioner degrades with the mesh disparity level.

\begin{table}[ht]
    \centering
    \begin{tabular}{c|lr|lr|lr|lr|lr|lr|lr}
	\toprule
	& \multicolumn{2}{c|}{$p=2$}
	& \multicolumn{2}{c|}{$p=3$}
	& \multicolumn{2}{c|}{$p=4$}
        & \multicolumn{2}{c|}{$p=5$}
        & \multicolumn{2}{c|}{$p=6$}
        & \multicolumn{2}{c|}{$p=7$}
	& \multicolumn{2}{c }{$p=8$}\\
	$d$
	& it & $\kappa$
	& it & $\kappa$
	& it & $\kappa$
        & it & $\kappa$
        & it & $\kappa$
        & it & $\kappa$
	& it & $\kappa$\\
	\midrule
        $1$ & 14 & 8.35 & 15 & 11.36 & 16 & 14.14 & 17 & 16.67 & 17 & 18.96 & 18 & 21.05 & 19 & 22.98\\
        $2$ & 14 & 8.53 & 16 & 11.57 & 16 & 14.24 & 17 & 16.79 & 18 & 19.09 & 20 & 21.19 & 20 & 23.13\\
        $3$ & 14 & 8.68 & 16 & 11.75 & 16 & 14.41 & 17 & 16.90 & 18 & 19.21 & 20 & 21.32 & 20 & 23.26\\
        $4$ & 14 & 8.82 & 16 & 11.92 & 16 & 14.60 & 17 & 17.01 & 18 & 19.33 & 20 & 21.45 & 20 & 23.39\\
        $5$ & 14 & 8.94 & 16 & 12.07 & 16 & 14.78 & 17 & 17.18 & 18 & 19.44 & 20 & 21.57 & 20 & 23.52\\
        $6$ & 14 & 9.06 & 16 & 12.22 & 16 & 14.95 & 17 & 17.36 & \multicolumn{2}{c|}{OoM} & \multicolumn{2}{c|}{OoM} & \multicolumn{2}{c}{OoM}\\
	\bottomrule
    \end{tabular}
    \caption{Iteration counts (it) and condition numbers $\kappa$ for the checkerboard mesh for varying disparity mesh level $d$ with scaled Dirichlet preconditioner with selection scaling.}
    \label{tab:checker_dis:selection}
\end{table}

\begin{table}[ht]
    \centering
    \begin{tabular}{c|lr|lr|lr|lr|lr|lr|lr}
	\toprule
	& \multicolumn{2}{c|}{$p=2$}
	& \multicolumn{2}{c|}{$p=3$}
	& \multicolumn{2}{c|}{$p=4$}
        & \multicolumn{2}{c|}{$p=5$}
        & \multicolumn{2}{c|}{$p=6$}
        & \multicolumn{2}{c|}{$p=7$}
	& \multicolumn{2}{c }{$p=8$}\\
	$d$
	& it & $\kappa$
	& it & $\kappa$
	& it & $\kappa$
        & it & $\kappa$
        & it & $\kappa$
        & it & $\kappa$
	& it & $\kappa$\\
	\midrule
        $1$ & 14 & 7.39 & 16 & 10.83 & 17 & 13.78 & 18 & 16.47 & 18 & 18.75 & 19 & 20.96 & 20 & 22.96\\
        $2$ & 15 & 9.50 & 17 & 13.09 & 18 & 16.25 & 19 & 19.07 & 20 & 21.62 & 20 & 23.76 & 21 & 25.84\\
        $3$ & 17 & 12.30 & 19 & 16.01 & 20 & 19.32 & 20 & 22.13 & 21 & 24.86 & 21 & 27.05 & 22 & 29.21\\
        $4$ & 18 & 15.63 & 20 & 19.52 & 21 & 22.88 & 21 & 25.72 & 22 & 28.50 & 23 & 31.03 & 24 & 33.14\\
        $5$ & 21 & 20.01 & 22 & 23.69 & 22 & 26.98 & 24 & 30.22 & 24 & 33.01 & 25 & 35.53 & 25 & 37.66\\
        $6$ & 22 & 25.00 & 23 & 28.28 & 25 & 32.02 & 25 & 35.07 & \multicolumn{2}{c|}{OoM} & \multicolumn{2}{c|}{OoM} & \multicolumn{2}{c}{OoM}\\
	\bottomrule
    \end{tabular}
    \caption{Iteration counts (it) and condition numbers $\kappa$ for the checkerboard mesh for varying disparity mesh level $d$ with deluxe type preconditioning.}
    \label{tab:checker_dis:edge}
\end{table}

\subsection{Adaptive refinement towards corner singularity}

In this example, we test the IETI-DP method with an adaptive refinement scheme with repeated patch-splitting as proposed in~\cite{Tyoler2025}. This leads to multi-patch configurations with emerging T-junctions. Again, the computational domain $\Omega$ is a quarter annulus. This time, it is composed of 4 patches. We again solve~\eqref{eq:nu poisson}, where we impose a diffusion coefficient of $\nu_{\text{orange}} = 1000$ for the lower left patch, whereas we leave $\nu_{\text{white}}=1$ for the remaining patches. This creates a singularity in the solution of the problem. To increase local approximation power, we use a residual-based a posteriori error estimator 
\begin{equation*}\label{res:err:est}
		\eta_k^2\coloneqq h_k^2 \|f+\Div(\nu\nabla u_h)\|_{L^2(\Omega_k)}^2
		+ \sum_{\ell \in N_\Gamma(k)} \frac{h_k}{2} \| \llbracket \nu \nabla u_h \rrbracket_\mathbf{n}\|_{L^2(\partial\Omega_k \cap \partial\Omega_\ell)}^2,
\end{equation*}
and use D\"orfler marking with $\theta=0.8$ to mark the specific patches for refinement where the estimator is too large, i.e., we choose a set of patch indices $\Mc \subset \{1,\dots,K\}$ for refinement such that $$\sum_{k \in \Mc} \eta_k^2 > \theta \sum_{k=1}^K \eta_k^2.$$

In Figure \ref{fig:ex2:mesh}, we present the initial multi-patch configuration and the adaptively refined mesh after sufficiently many refinement steps, as well as the computed solution $u$. Note that only the patches are visualized, the inner knot meshes are omitted. If we use the scaled Dirichlet preconditioner with selection scaling, we need to be careful since the condition number bound is not robust with respect to jumps in the diffusion coefficient. The following rule ensures that the solver does not degrade the condition number:
For every interface $\Gamma^{(k,\ell)}$ between $\Omega_k$ and $\Omega_\ell$ it holds that
\begin{equation}\label{eq:consistency}
        V_h^{(k)}|_{\Gamma^{(k,\ell)}} \subset V_h^{(\ell)}|_{\Gamma^{(k,\ell)}} \Rightarrow \nu_k \geq \nu_\ell.
\end{equation}

If we arrive at a mesh that does not satisfy \eqref{eq:consistency} or Assumption \ref{ass:primadmissible}, we split the respective patches until the assumptions are again satisfied (consistency splitting). Numerical evidence suggests that the number of additional patches obtained due to the consistency splitting is negligible. For the following results, we use a smaller tolerance and stop the iteration whenever the residual was reduced by a factor of $10^{-10}$. 

In Table \ref{tab:corner:selection:+consistency}, we present the results for the scaled Dirichlet preconditioner with selection scaling and consistency splitting. In Table \ref{tab:corner:selection:-consistency}, we present what happens if consistency splitting is not performed. One can observe that the iteration becomes somewhat less stable throughout the refinement process, since selection scaling sometimes uses the side with the smaller coefficient (the consistency condition~\eqref{eq:consistency} is not satisfied). Finally, in Table \ref{tab:corner:edge}, we present the results for the deluxe type preconditioner \eqref{eq:edgeprec}, where we again do not perform consistency splitting.

\newcommand{\s}{3.85}
\begin{figure}[ht]
\begin{minipage}[t]{0.3\textwidth}\centering
\includegraphics[height=\s cm]{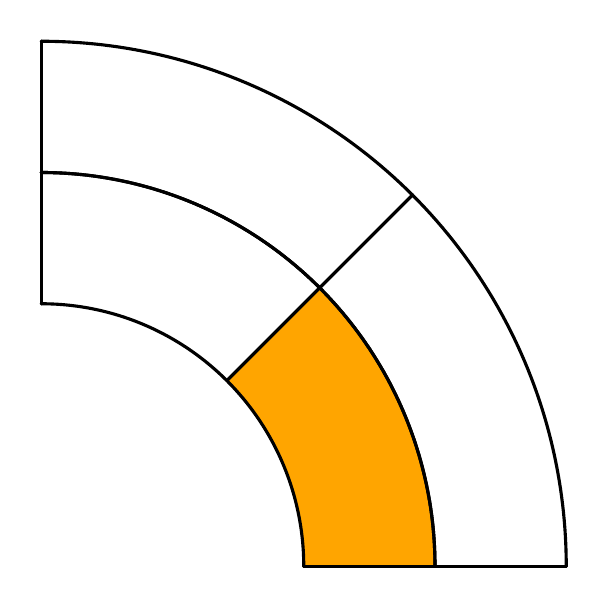}\\
\subcaption{Initial mesh}
\end{minipage}
\hfill
\begin{minipage}[t]{0.3\textwidth}\centering
\includegraphics[height=\s cm]{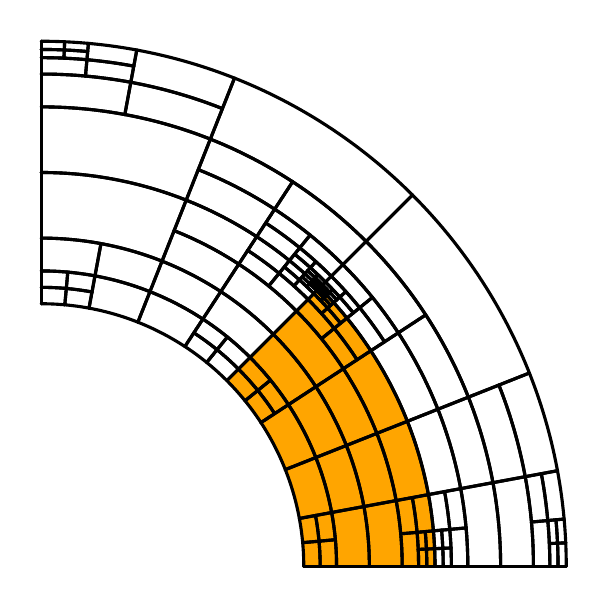}\\
\subcaption{Adaptive mesh}
\end{minipage}
\hfill
\begin{minipage}[t]{0.3\textwidth}\centering
\includegraphics[height=\s cm]{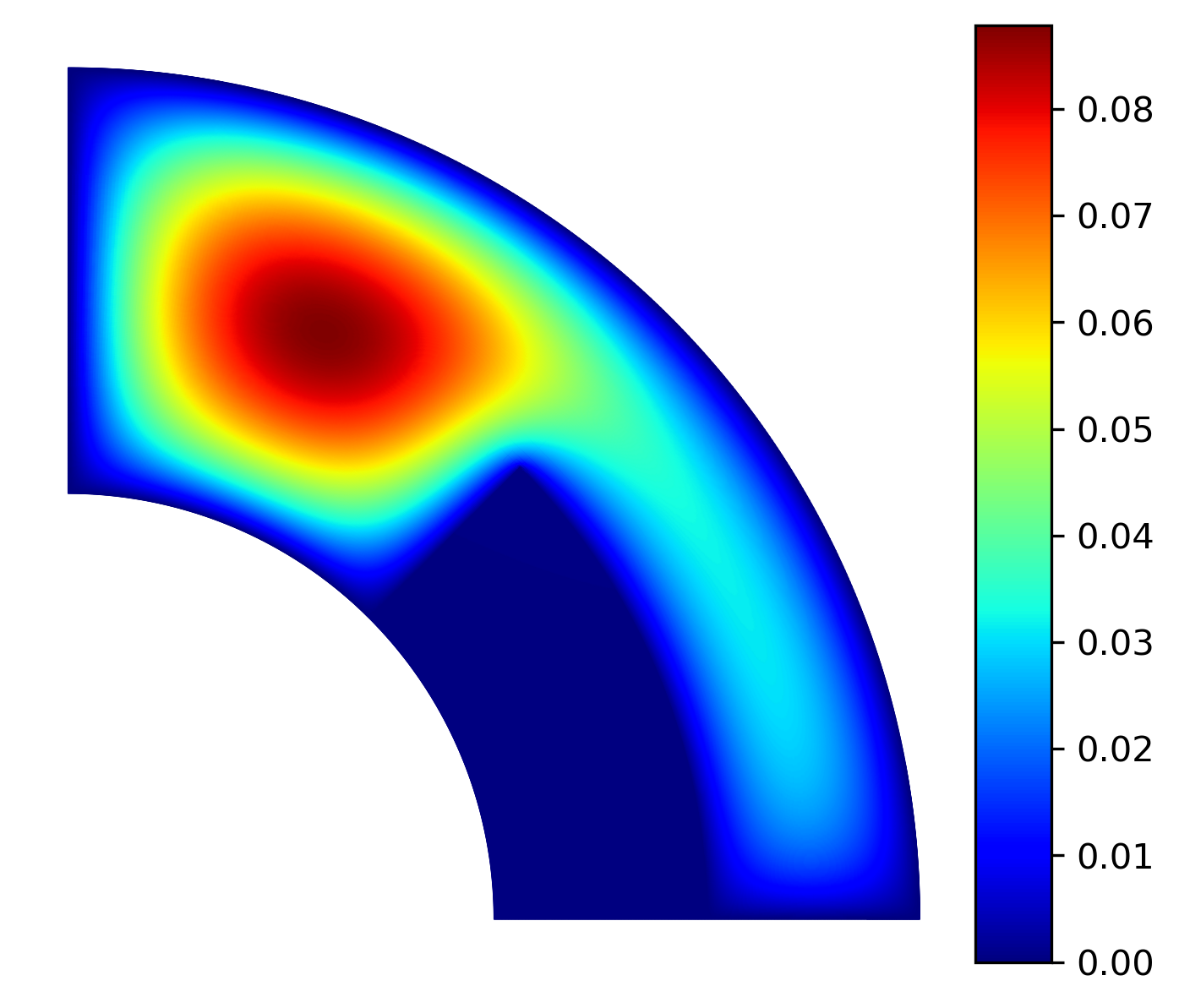}\\
\subcaption{Solution $u$}
\end{minipage}
\caption{\label{fig:ex2:mesh} Adaptive refinement form problem with corner singularity.}
\end{figure}

\begin{table}[ht]
    \centering
    \scriptsize 
    \setlength{\tabcolsep}{4.5pt}
    \begin{tabular}{c|lcr|lcr|lcr|lcr|lcr|lcr}
	\toprule
	& \multicolumn{3}{c|}{$p=2$}
	& \multicolumn{3}{c|}{$p=3$}
	& \multicolumn{3}{c|}{$p=4$}
    & \multicolumn{3}{c|}{$p=5$}
    & \multicolumn{3}{c|}{$p=6$}
    & \multicolumn{3}{c}{$p=7$}\\
	$r$
	& K & it & $\kappa$
	& K & it & $\kappa$
	& K & it & $\kappa$
    & K & it & $\kappa$
    & K & it & $\kappa$
	& K & it & $\kappa$\\
	\midrule
        1 & 4 & 7 & 2.72 & 4 & 8 & 3.24 & 4 & 8 & 3.70 & 4 & 9 & 4.10 & 4 & 9 & 4.47 & 4 & 9 & 4.80\\
        2 & 13 & 18 & 7.25 & 13 & 20 & 9.73 & 13 & 23 & 12.05 & 16 & 26 & 15.07 & 16 & 28 & 17.05 & 13 & 25 & 14.61\\
        3 & 28 & 26 & 7.81 & 25 & 28 & 10.76 & 31 & 33 & 12.70 & 28 & 35 & 16.39 & 28 & 37 & 18.59 & 22 & 36 & 17.95\\
        4 & 55 & 30 & 8.77 & 46 & 34 & 12.02 & 55 & 38 & 13.58 & 43 & 38 & 15.27 & 40 & 42 & 18.89 & 34 & 41 & 18.93\\
        5 & 79 & 30 & 8.89 & 67 & 34 & 12.37 & 79 & 39 & 13.80 & 61 & 41 & 15.64 & 64 & 43 & 18.94 & 46 & 44 & 19.19\\
        6 & 130 & 32 & 9.19 & 97 & 36 & 12.44 & 103 & 39 & 14.05 & 79 & 41 & 15.87 & 85 & 44 & 17.41 & 64 & 43 & 19.23\\
        7 & 217 & 33 & 9.40 & 130 & 36 & 12.49 & 127 & 40 & 14.26 & 106 & 42 & 16.83 & 103 & 44 & 17.46 & 88 & 45 & 19.20\\
        8 & 394 & 33 & 9.92 & 178 & 37 & 12.62 & 154 & 40 & 14.78 & 130 & 43 & 16.91 & 130 & 44 & 17.60 & 112 & 46 & 19.90\\
        9 & 601 & 34 & 10.42 & 253 & 38 & 12.67 & 181 & 40 & 14.86 & 160 & 43 & 16.98 & 148 & 45 & 17.59 & 136 & 47 & 19.94\\
        10 & 1186 & 35 & 10.45 & 388 & 39 & 13.01 & 238 & 41 & 15.00 & 205 & 44 & 17.08 & 178 & 45 & 17.59 & 163 & 47 & 20.79\\
        11 & 1921 & 35 & 11.12 & 670 & 40 & 13.80 & 385 & 43 & 16.26 & 316 & 46 & 18.74 & 229 & 46 & 19.06 & 208 & 48 & 20.94\\
        12 & & & & 1333 & 41 & 14.77 & 745 & 46 & 17.96 & 553 & 47 & 19.86 & 337 & 48 & 19.51 & 328 & 49 & 20.72\\
        13 &  & & & & & & & & & 967 & 48 & 19.47 & 520 & 50 & 21.85 & 559 & 51 & 23.24\\
	\bottomrule
    \end{tabular}
    \caption{\label{tab:corner:selection:+consistency} Iteration counts (it) and condition numbers $\kappa$ for the L-shape domain after $r$ adaptive refinements with scaled Dirichlet preconditioner with selection scaling and additional consistency refinement.}
\end{table}

\begin{table}[ht]
    \centering
    \scriptsize 
    \setlength{\tabcolsep}{3.8pt}
    \begin{tabular}{c|lcr|lcr|lcr|lcr|lcr|lcr}
	\toprule
	& \multicolumn{3}{c|}{$p=2$}
	& \multicolumn{3}{c|}{$p=3$}
	& \multicolumn{3}{c|}{$p=4$}
    & \multicolumn{3}{c|}{$p=5$}
    & \multicolumn{3}{c|}{$p=6$}
    & \multicolumn{3}{c}{$p=7$}\\
	$r$
	& K & it & $\kappa$
	& K & it & $\kappa$
	& K & it & $\kappa$
    & K & it & $\kappa$
    & K & it & $\kappa$
	& K & it & $\kappa$\\
	\midrule
        1 & 4 & 7 & 2.72 & 4 & 8 & 3.24 & 4 & 8 & 3.70 & 4 & 9 & 4.10 & 4 & 9 & 4.47 & 4 & 9 & 4.80\\
        2 & 13 & 18 & 7.25 & 13 & 20 & 9.73 & 13 & 23 & 12.05 & 13 & 69 & 2.9e4 & 13 & 84 & 3.2e4 & 10 & 74 & 3.5e4\\
        3 & 28 & 26 & 7.81 & 25 & 28 & 10.76 & 31 & 33 & 12.70 & 22 & 70 & 2.3e4 & 19 & 79 & 3.2e4 & 22 & 133 & 4.1e4\\
        4 & 55 & 30 & 8.77 & 46 & 34 & 12.02 & 52 & 77 & 2.6e4 & 31 & 36 & 14.97 & 28 & 101 & 3.2e4 & 28 & 134 & 3.3e4\\
        5 & 79 & 30 & 8.89 & 64 & 60 & 1.6e4 & 70 & 113 & 2.9e4 & 43 & 122 & 3.3e4 & 34 & 150 & 3.7e4 & 40 & 152 & 3.3e4\\
        6 & 130 & 32 & 9.19 & 85 & 35 & 11.38 & 91 & 89 & 1.9e4 & 61 & 139 & 3.2e4 & 43 & 144 & 2.9e4 & 58 & 164 & 3.3e4\\
        7 & 214 & 48 & 1.2e4 & 106 & 36 & 11.71 & 112 & 39 & 14.05 & 85 & 147 & 3.2e4 & 58 & 138 & 2.7e4 & 79 & 186 & 3.0e4\\
        8 & 391 & 47 & 1.8e4 & 142 & 93 & 2.2e4 & 127 & 39 & 14.26 & 109 & 151 & 2.6e4 & 73 & 152 & 2.7e4 & 100 & 176 & 2.9e4\\
        9 & 595 & 60 & 1.8e4 & 199 & 100 & 2.2e4 & 151 & 40 & 14.42 & 133 & 151 & 2.6e4 & 97 & 149 & 2.7e4 & 118 & 183 & 2.9e4\\
        10 & 1147 & 101 & 1.9e4 & 286 & 104 & 2.2e4 & 181 & 40 & 14.46 & 157 & 154 & 2.6e4 & 115 & 154 & 2.7e4 & 136 & 189 & 2.9e4\\
        11 & 1882 & 113 & 2.0e4 & 457 & 131 & 2.2e4 & 229 & 85 & 2.0e4 & 205 & 156 & 2.6e4 & 139 & 158 & 2.7e4 & 166 & 188 & 2.9e4\\
        12 & & & & 847 & 163 & 2.2e4 & 358 & 163 & 2.6e4 & 304 & 197 & 2.8e4 & 163 & 160 & 2.7e4 & 211 & 189 & 2.9e4\\
        13 & & & & 1669 & 201 & 2.2e4 & 709 & 215 & 2.8e4 & 562 & 242 & 3.0e4 & 208 & 162 & 2.7e4 & & & \\
	\bottomrule
    \end{tabular}
    \caption{\label{tab:corner:selection:-consistency}Iteration counts (it) and condition numbers $\kappa$ for the L-shape domain after $r$ adaptive refinements with scalaed Dirichlet preconditioner with selection scaling and without additional consistency refinement.}
\end{table}

\begin{table}[ht]
    \centering
    \scriptsize 
    \setlength{\tabcolsep}{4.5pt}
    \begin{tabular}{c|lcr|lcr|lcr|lcr|lcr|lcr}
	\toprule
	& \multicolumn{3}{c|}{$p=2$}
	& \multicolumn{3}{c|}{$p=3$}
	& \multicolumn{3}{c|}{$p=4$}
    & \multicolumn{3}{c|}{$p=5$}
    & \multicolumn{3}{c|}{$p=6$}
    & \multicolumn{3}{c}{$p=7$}\\
	$r$
	& K & it & $\kappa$
	& K & it & $\kappa$
    & K & it & $\kappa$
    & K & it & $\kappa$
    & K & it & $\kappa$
	& K & it & $\kappa$\\
	\midrule
        1 & 4 & 8 & 2.65 & 4 & 9 & 3.74 & 4 & 10 & 4.69 & 4 & 11 & 5.52 & 4 & 12 & 6.28 & 4 & 12 & 6.96\\
        2 & 13 & 19 & 6.26 & 13 & 22 & 9.25 & 13 & 25 & 12.11 & 13 & 25 & 13.15 & 13 & 27 & 15.19 & 10 & 26 & 15.12\\
        3 & 28 & 25 & 7.06 & 25 & 28 & 10.26 & 31 & 34 & 13.29 & 22 & 32 & 15.31 & 19 & 33 & 16.31 & 22 & 38 & 18.66\\
        4 & 55 & 27 & 6.96 & 46 & 32 & 10.64 & 52 & 36 & 13.22 & 31 & 38 & 16.06 & 28 & 39 & 16.95 & 28 & 40 & 18.93\\
        5 & 79 & 27 & 6.95 & 64 & 32 & 10.25 & 70 & 36 & 13.38 & 43 & 39 & 15.42 & 34 & 41 & 16.86 & 40 & 44 & 18.95\\
        6 & 130 & 27 & 6.87 & 85 & 33 & 10.34 & 91 & 37 & 12.54 & 61 & 40 & 15.41 & 43 & 42 & 16.95 & 58 & 45 & 18.97\\
        7 & 214 & 28 & 6.99 & 106 & 33 & 10.32 & 112 & 37 & 13.31 & 85 & 41 & 15.42 & 58 & 43 & 17.05 & 79 & 46 & 19.01\\
        8 & 391 & 28 & 7.25 & 142 & 33 & 9.87 & 127 & 38 & 13.32 & 109 & 40 & 15.32 & 73 & 43 & 17.29 & 100 & 45 & 20.12\\
        9 & 595 & 29 & 7.28 & 199 & 34 & 10.45 & 151 & 38 & 13.25 & 133 & 41 & 15.34 & 97 & 43 & 17.51 & 118 & 46 & 20.12\\
        10 & 1147 & 29 & 7.40 & 286 & 34 & 10.68 & 181 & 38 & 13.25 & 157 & 41 & 15.33 & 115 & 43 & 17.49 & 136 & 47 & 20.06\\
        11 & 1882 & 29 & 7.35 & 457 & 35 & 10.98 & 229 & 39 & 13.32 & 205 & 41 & 15.49 & 139 & 44 & 17.35 & 166 & 46 & 19.56\\
        12 & & & & 847 & 36 & 11.15 & 358 & 39 & 13.33 & 304 & 42 & 16.20 & 163 & 44 & 17.38 & 211 & 46 & 20.17\\
        13 & & & & 1669 & 36 & 11.21 & 709 & 40 & 14.15 & 562 & 44 & 16.87 & 208 & 44 & 17.92 & & &\\
	\bottomrule
    \end{tabular}
    \caption{\label{tab:corner:edge}Iteration counts (it) and condition numbers $\kappa$ for the L-shape domain after $r$ adaptive refinements with deluxe type preconditioner.}
\end{table}

\subsection{Simulation of inductor machine with adaptive refinement}

As last example, we apply the IETI-DP method to the simulation of a gapped inductor, where Maxwell's equations model the electromagnetic effects. We consider linear magneto-statics, so for given permeability $\mu$ and current density $J$, the magnetic field density $H$ and magnetic flux density $B$ are given by
\begin{align*}
    \Curl\,H = J, \quad
    \Div\,B = 0, \quad
    B = \mu H.
\end{align*}
We model the cross-section of the gapped inductor, consisting of an E-shaped iron core, 2 slots for the copper windings, and an iron piece distanced from the E-shaped part by a small air gap. The entire computational domain encompasses the machine and an empty box of air surrounding the machine. The 3-dimensional domain $\Omega^* = \Omega \times [0,L]$ is reduced to its 2-dimensional cross-section $\Omega$, see Figure~\ref{fig:ex2:mesh}. Since the domain is simply-connected, the condition $\Div B = 0$ implies that there exists a vector potential $A$ such that $B = \Curl A$. Since current is only flowing in transversal direction of the cross-section, i.e., $J = (0,0,j(x,y))^\top$, we further assume the magnetic flux and field strength to have no component in this direction and to be independent transversal position. In total, we obtain $H=(H_x(x,y),H_y(x,y),0)^\top$, $B=(B_x(x,y),B_y(x,y),0)^\top$ and $A=(0,0,u(x,y))^\top$ which reduces the whole problem to a scalar partial differential equation on the cross-section of the inductor. Here, Maxwell's equations simplify to: Find $u \in H_0^1(\Omega)$ such that
\begin{alignat*}{2}
    \int_\Omega \mu^{-1} \nabla u \cdot \nabla v \diff x
    = \int_\Omega j v \diff x \quad \text{for all } v \in H_0^1(\Omega).
\end{alignat*}
The current density $j$ is nonzero only in the copper windings, where we choose $1$ Ampere per winding. We have $66$ windings and a total of $66$ Ampere, which we divide by the area of the winding slots ($0.008\,m^2$) to get the current density $j= 66/0.008$ Ampere/$m^2$. The direction of the current in the left slot is positive and in the right slot negative. The permeability is given by $\mu_{\text{Air}} = \mu_{\text{Copper}} = 4\pi \cdot 10^{-7}$ and $\mu_{\text{Iron}} = 204\pi \cdot 10^{-5}$.

In Figure \ref{fig:Ind:mesh}, we show the initial multi-patch configuration consisting of $42$ patches and the adaptively refined mesh after sufficiently many refinement steps, as well as the magnitude of the computed magnetic flux density $|B| = |\mu^{-1}\nabla^\perp u|$ and the computed magnetic field strength $|H| = \mu^{-1}|B|$. Observe that $H$ is large inside the air gap and has singularities radiating around the corners near the air gap. Further, we see how the error estimator detects the singularities and refines the specific corners.  

In Table \ref{tab:Ind:selection:+consistency}, we present the results for the scaled Dirichlet preconditioner with selection scaling, where consistency splits have been performed in order to satisfy Assumption \ref{ass:primadmissible} as well as the consistency condition \eqref{eq:consistency}. We can see that the condition number is slightly larger than for the quarter annulus example, since there is a large aspect ratio at the air gap, which degrades the constant $C_G$. This can be easily solved by splitting the patches in an anisotropic manner beforehand, which we omitted in order to start with a fully matching patch configuration.

\newcommand{\sss}{3.8}
\begin{figure}[htb]
\begin{minipage}[t]{0.47\textwidth}\centering
\includegraphics[height=\sss cm]{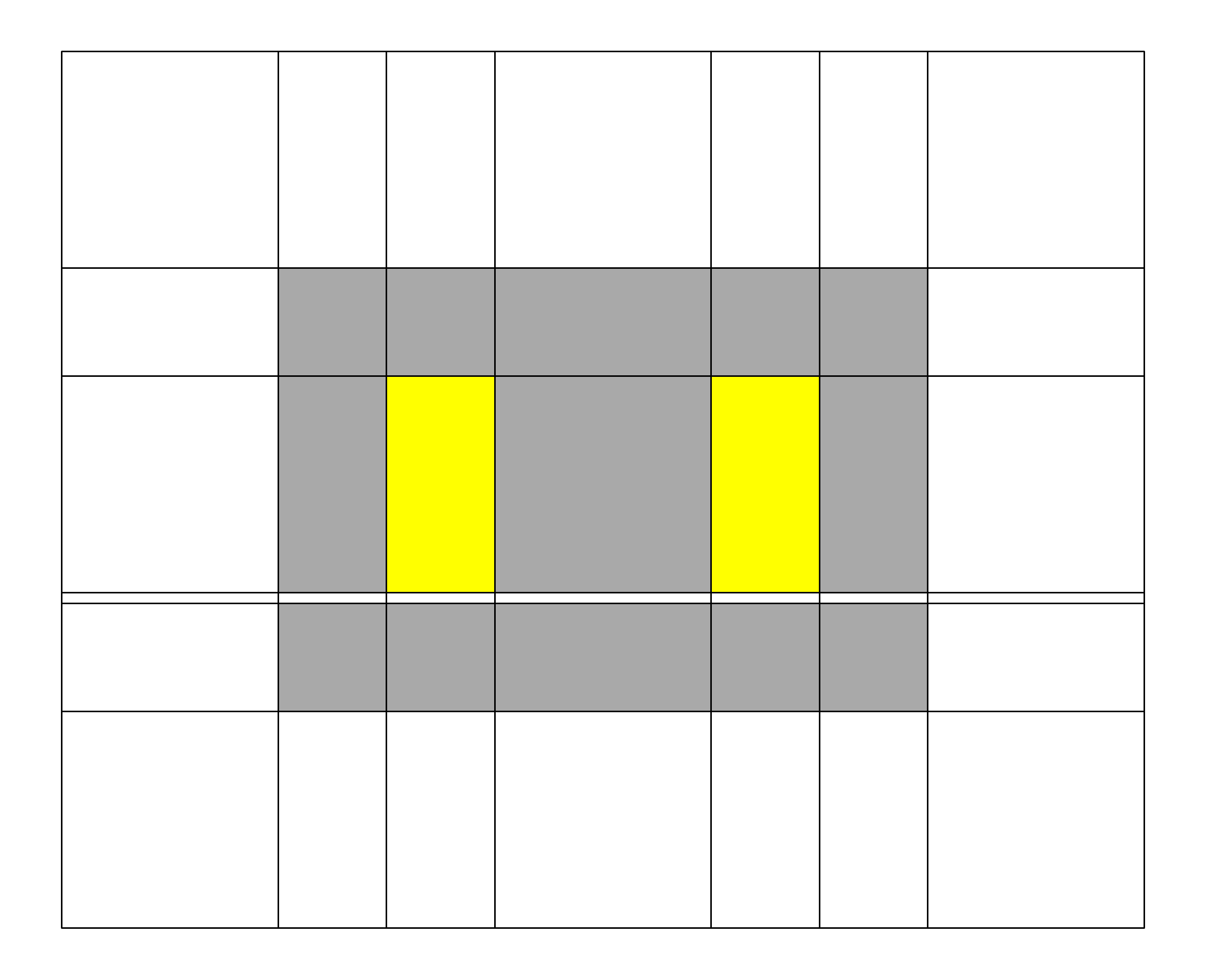}\\
\subcaption{Material distribution of the Inductor}
\end{minipage}
\hfill
\begin{minipage}[t]{0.47\textwidth}\centering
\includegraphics[height=\sss cm]{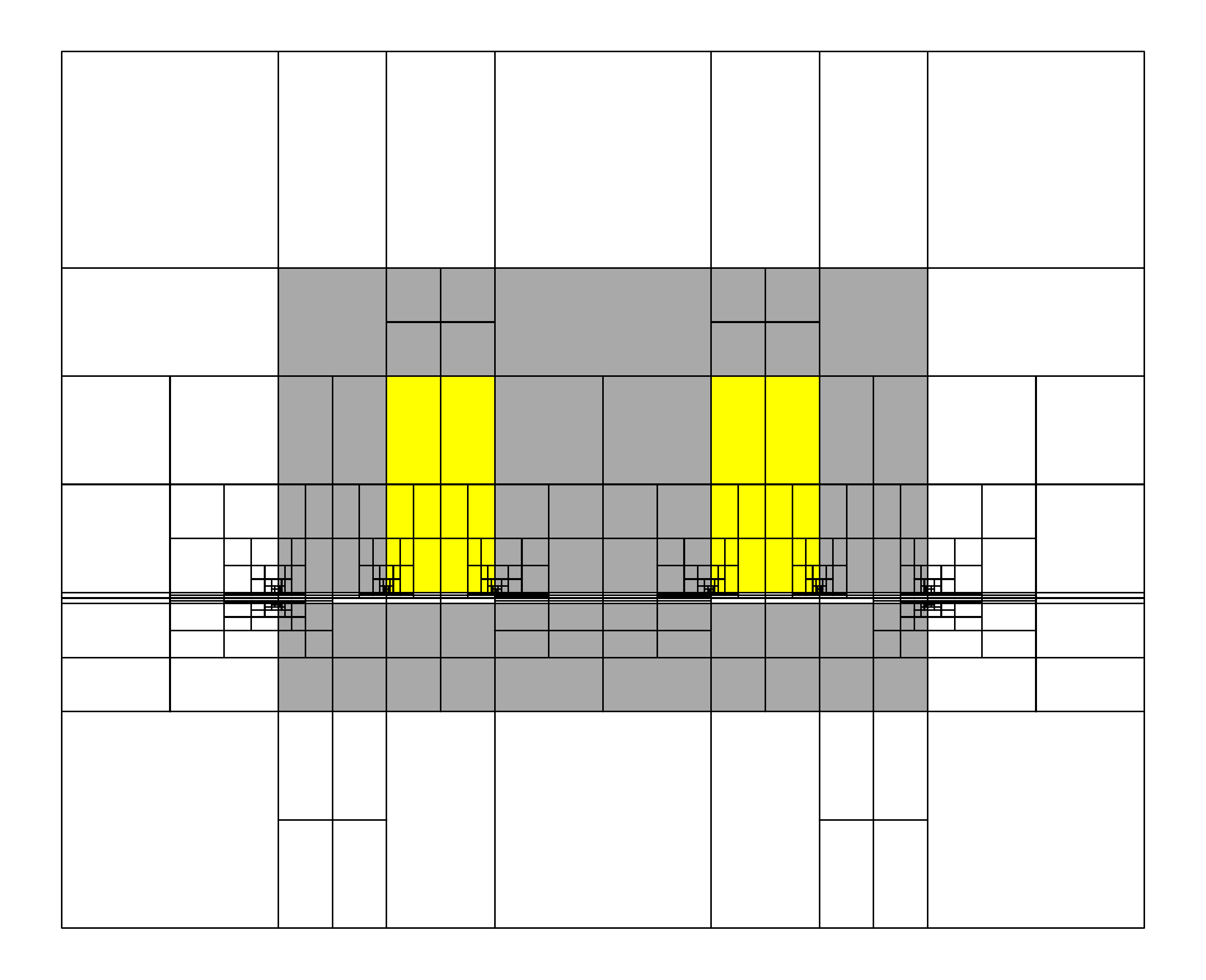}\\
\subcaption{Multi-patch domain after adaptive refinement}
\end{minipage}
\vskip\baselineskip
\begin{minipage}[t]{0.47\textwidth}\centering
\includegraphics[height=\sss cm]{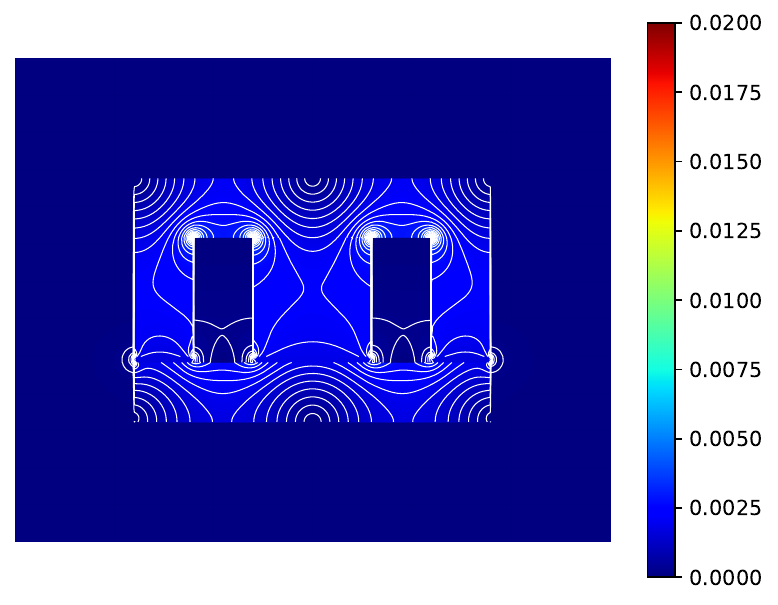}\\
\subcaption{The magnitude of the magnetic flux density B}
\end{minipage}
\hfill
\begin{minipage}[t]{0.47\textwidth}\centering
\includegraphics[height=\sss cm]{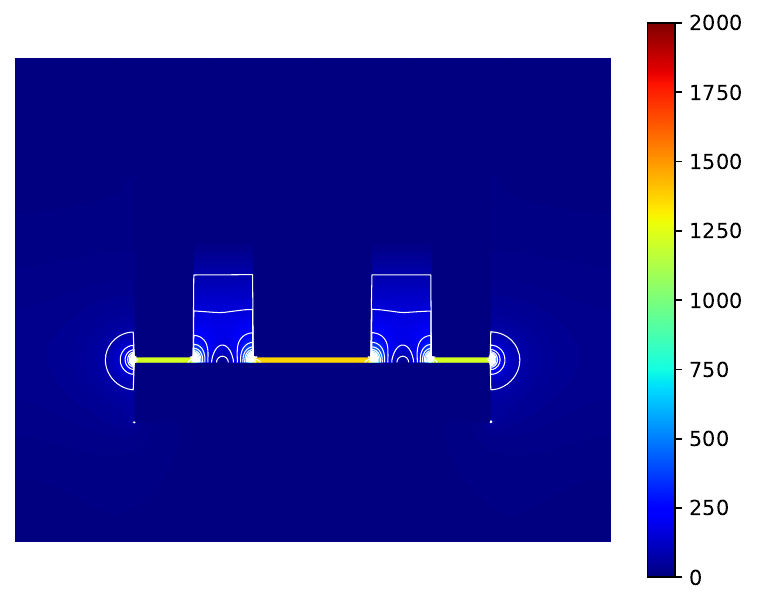}\\
\subcaption{The magnitude of the magnetic field H}
\end{minipage}
\caption{\label{fig:Ind:mesh} Inductor machine layout and computed fields.}
\end{figure}

\begin{table}[ht]
    \centering
    \scriptsize 
    \setlength{\tabcolsep}{3.2pt}
    \begin{tabular}{c|lcr|lcr|lcr|lcr|lcr|lcr}
	\toprule
	& \multicolumn{3}{c|}{$p=2$}
	& \multicolumn{3}{c|}{$p=3$}
	& \multicolumn{3}{c|}{$p=4$}
    & \multicolumn{3}{c|}{$p=5$}
    & \multicolumn{3}{c|}{$p=6$}
    & \multicolumn{3}{c}{$p=7$}\\
	$r$
	& K & it & $\kappa$
	& K & it & $\kappa$
	& K & it & $\kappa$
    & K & it & $\kappa$
    & K & it & $\kappa$
	& K & it & $\kappa$\\
	\midrule
        1 & 42 & 30 & 57.25 & 42 & 35 & 59.11 & 42 & 37 & 60.87 & 42 & 38 & 62.47 & 42 & 39 & 63.87 & 42 & 40 & 65.10\\
        2 & 99 & 60 & 95.02 & 99 & 64 & 96.28 & 87 & 66 & 107.29 & 75 & 67 & 111.25 & 75 & 70 & 115.33 & 75 & 72 & 118.71\\
        3 & 183 & 98 & 143.57 & 177 & 94 & 142.13 & 153 & 98 & 150.10 & 135 & 99 & 154.74 & 135 & 104 & 160.20 & 135 & 107 & 163.47\\
        4 & 261 & 106 & 144.42 & 255 & 98 & 139.38 & 231 & 104 & 145.81 & 219 & 111 & 148.10 & 213 & 115 & 159.93 & 207 & 116 & 162.85\\
        5 & 363 & 119 & 153.05 & 345 & 112 & 143.75 & 309 & 112 & 139.70 & 297 & 116 & 150.34 & 291 & 118 & 161.85 & 279 & 120 & 162.46\\
        6 & 519 & 126 & 164.47 & 441 & 116 & 143.06 & 387 & 116 & 139.52 & 381 & 119 & 149.42 & 369 & 121 & 157.57 & 363 & 122 & 160.87\\
        7 & 681 & 136 & 166.88 & 561 & 124 & 150.44 & 477 & 119 & 139.44 & 459 & 123 & 153.89 & 453 & 123 & 157.45 & 441 & 124 & 160.83\\
        8 & 915 & 142 & 175.82 & 753 & 137 & 171.50 & 585 & 127 & 155.15 & 549 & 126 & 144.58 & 531 & 127 & 157.45 & 531 & 128 & 160.79\\
        9 & & & & 1011 & 145 & 174.05 & 693 & 131 & 155.69 & 645 & 128 & 144.28 & 627 & 128 & 157.38 & 621 & 130 & 160.73\\
        10& & & & & & & 849 & 139 & 173.90 & 765 & 134 & 160.24 & 723 & 129 & 157.36 & 717 & 131 & 160.66\\
        11& & & &  &  &  &  &  &  & 915 & 141 & 177.09 & 837 & 129 & 157.35 & 819 & 132 & 160.64\\
	\bottomrule
    \end{tabular}
    \caption{\label{tab:Ind:selection:+consistency} Iteration counts (it) and condition numbers $\kappa$ for the inductor machine after $r$ adaptive refinements with selection scaling and additional consistency refinement.}
\end{table}

\subsection{Conclusions}

We have applied the IETI-DP method to discretizations that are potentially not fully matching at the interfaces, but with nested trace spaces. This assumption is naturally satisfied by the adaptive refinement strategy proposed in~\cite{Tyoler2025}. The function values at the T-junctions are chosen as additional primal degrees of freedom in addition to the usual corner vertices. This choice offers two advantages: It decouples the constraints handled via the Lagrange multipliers, and it requires only a single primal degree of freedom per T-junction, thereby ensuring a $p$-robust number of primal degrees of freedom. With Theorem~\ref{thrm:final}, we have established a condition number bound of order $p(1 + \log p + \log \tfrac{H}{h})^2$, matching the optimal estimates known for classical FETI-DP and IETI-DP methods.

If the diffusion coefficient are not constant, Theorem~\ref{thrm:final} indicates that the condition number might also depend on the jumps of those coefficients, which can also be seen in the numerical experiments. We have outlined two simple remedies: first, applying a deluxe-type preconditioner to the affected interfaces, which is more expensive, but restores robustness. The second option is to apply additional refinements. The numerical examples indicate that this strategy only introduces a few new patches, which does not affect optimal convergence rates. 

Our numerical experiments confirm both the predicted degradation in the presence of diffusion jumps and the effectiveness of the proposed remedies. In particular, additional adaptive refinement and the use of deluxe-type preconditioner successfully mitigate the dependence on coefficient jumps, thereby extending the robustness and applicability of IETI-DP to adaptively refined multi-patch discretizations.

\subsection*{Acknowledgements}
This research was funded in whole or in part by the Austrian Science Fund (FWF):10.55776/P33956. 

\bibliographystyle{plain}


\begin{thebibliography}{10}

\bibitem{Adams2003}
Robert~A. Adams and John J.~F. Fournier.
\newblock {\em Sobolev spaces}.
\newblock Pure and applied mathematics. Elsevier Academic Press, 2nd edition, 2003.

\bibitem{Brivadis2015}
Ericka Brivadis, Annalisa Buffa, Barbara Wohlmuth, and Linus Wunderlich.
\newblock Isogeometric mortar methods.
\newblock {\em Computer Methods in Applied Mechanics and Engineering}, 284:292--319, 2015.
\newblock Isogeometric Analysis Special Issue.

\bibitem{Pavarino2013}
Louren{\c{c}}o~Beir{\~a}o da~Veiga, David Cho, Luca~F. Pavarino, and Simone Scacchi.
\newblock {BDDC} preconditioners for {I}sogeometric {A}nalysis.
\newblock {\em Mathematical Models and Methods in Applied Sciences}, 23(06):1099--1142, 2013.

\bibitem{Widlund2014}
Louren{\c{c}}o~Beir{\~a}o da~Veiga, Luca~F. Pavarino, Simone Scacchi, Olof~B. Widlund, and Stefano Zampini.
\newblock Isogeometric {BDDC} preconditioners with deluxe scaling.
\newblock {\em SIAM Journal on Scientific Computing}, 36(3):A1118--A1139, 2014.

\bibitem{deBoor1972}
Carl de~Boor.
\newblock {On calculating with B-splines}.
\newblock {\em Journal of Approximation Theory}, 6:50--62, 7 1972.

\bibitem{deBoor1978}
Carl de~Boor.
\newblock {\em A Practical Guide to Splines}.
\newblock Applied Mathematical Sciences. Springer, New York, revised edition, 1978.

\bibitem{Dokken2013}
Tor Dokken, Tom Lyche, and Kjell~Fredrik Pettersen.
\newblock Polynomial splines over locally refined box-partitions.
\newblock {\em Computer Aided Geometric Design}, 30(3):331--356, 2013.

\bibitem{Dörfel2010}
Michael~R. Dörfel, Bert Jüttler, and Bernd Simeon.
\newblock Adaptive {I}sogeometric {A}nalysis by local h-refinement with {T}-splines.
\newblock {\em Computer Methods in Applied Mechanics and Engineering}, 199(5):264--275, 2010.
\newblock Computational Geometry and Analysis.

\bibitem{Doerfler1996}
Willy Dörfler.
\newblock A convergent adaptive algorithm for {P}oisson's equation.
\newblock {\em SIAM Journal on Numerical Analysis}, 33:1106--1124, 1996.

\bibitem{Farhat1991}
Charbel Farhat and Francois-Xavier Roux.
\newblock A method of finite element tearing and interconnecting and its parallel solution algorithm.
\newblock {\em International Journal for Numerical Methods in Engineering}, 32(6):1205--1227, 1991.

\bibitem{Giannelli2012}
Carlotta Giannelli, Bert Jüttler, and Hendrik Speleers.
\newblock {THB}-splines: The truncated basis for hierarchical splines.
\newblock {\em Computer Aided Geometric Design}, 29:485--498, 10 2012.

\bibitem{Hofer2017}
Christoph Hofer and Ulrich Langer.
\newblock {Dual-primal isogeometric tearing and interconnecting solvers for multipatch dG-IgA equations}.
\newblock {\em Computer Methods in Applied Mechanics and Engineering}, 316:2--21, 2017.

\bibitem{Hofreither2022}
Clemens Hofreither, Ludwig Mitter, and Hendrik Speleers.
\newblock Local multigrid solvers for adaptive {I}sogeometric {A}nalysis in hierarchical spline spaces.
\newblock {\em IMA Journal of Numerical Analysis}, 42:2429--2458, 2022.

\bibitem{Hughes2005}
Thomas J.~R. Hughes, John~A. Cottrell, and Yuri Bazilevs.
\newblock Isogeometric analysis: {CAD}, finite elements, {NURBS}, exact geometry and mesh refinement.
\newblock {\em Computer Methods in Applied Mechanics and Engineering}, 194:4135--4195, 2005.

\bibitem{Dokken2014}
Kjetil~A. Johannessen, Trond Kvamsdal, and Tor Dokken.
\newblock Isogeometric analysis using {LR B}-splines.
\newblock {\em Computer Methods in Applied Mechanics and Engineering}, 269:471--514, 2014.

\bibitem{Kleiss2012}
Stefan Kleiss, Clemens Pechstein, Bert Jüttler, and Satyendra Tomar.
\newblock {IETI} – isogeometric tearing and interconnecting.
\newblock {\em Computer Methods in Applied Mechanics and Engineering}, 247-248:201--215, 2012.

\bibitem{Kraft1997}
Rainer Kraft.
\newblock Adaptive and linearly independent multilevel {B}-splines.
\newblock In A.~Le~Méhauté, C.~Rabut, and L.~L. Schumaker, editors, {\em Surface Fitting and Multiresolution Methods}, pages 209 -- 218. Vanderbilt University Press, Nashville, 1997.

\bibitem{LangerToulopoulos2015}
Ulrich Langer and Ioannis Toulopoulos.
\newblock {Analysis of multipatch discontinuous Galerkin IgA approximations to elliptic boundary value problems}.
\newblock {\em Computing and Visualization in Science}, 17:217–233, 2015.

\bibitem{Mandel2005}
Jan Mandel, Clark~R. Dohrmann, and Radek Tezaur.
\newblock {An algebraic theory for primal and dual substructuring methods by constraints}.
\newblock {\em Applied Numerical Mathematics}, 54(2):167--193, 2005.
\newblock 6th IMACS.

\bibitem{Praetorius2020}
Carl-Martin Pfeiler and Dirk Praetorius.
\newblock Dörfler marking with minimal cardinality is a linear complexity problem.
\newblock {\em Mathematics of Computation}, 89:2735--2752, 2020.

\bibitem{Riviere2008}
Béatrice Rivière.
\newblock {\em Discontinuous Galerkin Methods for Solving Elliptic and Parabolic Equations}.
\newblock Society for Industrial and Applied Mathematics, 2008.

\bibitem{Schneckenleitner2020}
Rainer Schneckenleitner and Stefan Takacs.
\newblock {Condition number bounds for IETI-DP methods that are explicit in $p$ and $h$}.
\newblock {\em Mathematical Models and Methods in Applied Sciences}, 30(11):2067--2103, 2020.

\bibitem{Schneckenleitner2022}
Rainer Schneckenleitner and Stefan Takacs.
\newblock {IETI-DP} methods for discontinuous {Galerkin} multi-patch {Isogeometric Analysis} with {T}-junctions.
\newblock {\em Computer Methods in Applied Mechanics and Engineering}, 393:114694, 2022.

\bibitem{Schwab1998}
Christoph Schwab.
\newblock {\em $p$- and $hp$-Finite Element Methods: Theory and Applications in Solid and Fluid Mechanics}.
\newblock Clarendon Press Oxford, 1998.

\bibitem{Sederberg2012}
Michael~A. Scott, Xin Li, Thomas~W. Sederberg, and Thomas J.~R. Hughes.
\newblock Local refinement of analysis-suitable {T}-splines.
\newblock {\em Computer Methods in Applied Mechanics and Engineering}, 213-216:206--222, 2012.

\bibitem{Tyoler2025}
Stefan Takacs and Stefan Tyoler.
\newblock Multi-resolution {I}sogeometric {A}nalysis – efficient adaptivity utilizing the multi-patch structure.
\newblock {\em Computers \& Mathematics with Applications}, 179:103--125, 2025.

\bibitem{Toselli2005}
Andrea Toselli and Olof~B. Widlund.
\newblock {\em Domain Decomposition Methods - Algorithms and Theory}.
\newblock Springer Series in Computational Mathematics. Springer Berlin, Heidelberg, 1st edition edition, 2005.

\bibitem{Verfuerth2013}
Rüdiger Verfürth.
\newblock {\em A posteriori error estimation techniques for finite element methods}.
\newblock Oxford University Press, 2013.

\bibitem{Giannelli2011}
Anhvu Vuong, Carlotta Giannelli, Bert Jüttler, and Bernd Simeon.
\newblock A hierarchical approach to adaptive local refinement in {I}sogeometric {A}nalysis.
\newblock {\em Computer Methods in Applied Mechanics and Engineering}, 200(49):3554--3567, 2011.

\end{thebibliography}

\end{document}